\documentclass[11pt]{article}

\usepackage{amsmath, amsfonts, amssymb,latexsym,wasysym,graphicx}
\input epsf.sty

\newtheorem{Theorem}{Theorem}[section]
\newtheorem{Lemma}{Lemma}[section]
\newtheorem{Proposition}[Lemma]{Proposition}

\newtheorem{Definition}[Lemma]{Definition}

\newtheorem{Example}[Lemma]{Example}

\newcommand{\BEQ}{\begin{equation}}     
\newcommand{\BEA}{\begin{eqnarray}}
\newcommand{\BD}{\begin{displaymath}}
\newcommand{\EEQ}{\end{equation}}       
\newcommand{\EEA}{\end{eqnarray}}
\newcommand{\ED}{\end{displaymath}}

\newcommand{\del}{\delta}
\newcommand{\Del}{\Delta}
\newcommand{\eps}{\varepsilon}          

\newcommand{\supp}{{\mathrm{supp}}}
\newcommand{\R}{\mathbb{R}}
\newcommand{\C}{\mathbb{C}}

\newcommand{\N}{\mathbb{N}}

\newcommand{\Id}{{\mathrm{Id}}}
\newcommand{\Sk}{{\mathrm{Sk}}}
\newcommand{\SkI}{{\mathrm{Sk\,I}}}

\def\esper{{\mathbb{E}}}
\def\T{{\mathbb{T}}}
\def\Var{{\mathrm{Var}}}
\def\Cov{{\mathrm{Cov}}}

\def\sgn{{\mathrm{sgn}}}

\newcommand{\eop}{\hfill $\Box$}        

\newcommand{\II}{{\rm i}}               
\renewcommand{\Re}{{\rm Re\ }}          
\renewcommand{\Im}{{\rm Im\ }}          
\newcommand{\half}{{1\over 2}}          
\renewcommand{\vec}[1]{\boldsymbol{#1}} 

\catcode`\@=11
\def\numberbysection{\@addtoreset{equation}{section}
        \def\theequation{\thesection.\arabic{equation}}}
\numberbysection

\begin{document}

\vspace*{1.5cm}
\begin{center}
{\Large \bf A rough path over
multidimensional fractional Brownian motion with arbitrary Hurst index by
Fourier normal ordering
}
\end{center}

\vspace{2mm}
\begin{center}
{\bf  J\'er\'emie Unterberger}
\end{center}

\vspace{2mm}
\begin{quote}

\renewcommand{\baselinestretch}{1.0}
\footnotesize
{Fourier normal ordering \cite{Unt09bis} is a new algorithm to construct explicit rough paths
over arbitrary H\"older-continuous multidimensional paths. We apply in this article
the Fourier normal ordering ordering algorithm to the construction of 
 an explicit rough path over  multi-dimensional fractional
Brownian motion $B$  with arbitrary Hurst index $\alpha$ (in particular, for $\alpha\le 1/4$,
which was till now an open problem)
 by regularizing the iterated integrals of the analytic approximation of $B$
 defined in \cite{Unt08}.  The regularization procedure is applied to 'Fourier normal ordered'
iterated integrals obtained by permuting the order of integration so that innermost integrals
have highest Fourier modes. The algebraic properties of this rough path are  best understood 
 using two Hopf algebras: the Hopf algebra of decorated
rooted trees \cite{ConKre98} for the multiplicative or Chen property, and the
shuffle algebra for the geometric or shuffle property. The rough path lives in Gaussian
chaos of integer orders and is shown to have finite moments.

As well-known, the construction of a 
 rough path is the key to  
defining a stochastic calculus  and solve stochastic differential equations driven by
$B$.

 The article \cite{Unt09ter} gives a quick overview of the
method.

 }
\end{quote}

\vspace{4mm}


\tableofcontents
\newpage


\section{Introduction}


The (two-sided) fractional Brownian motion $t\to B_t$, $t\in\R$  (fBm for short) with Hurst exponent $\alpha$, $\alpha\in(0,1)$, defined as the centered Gaussian process with covariance
\BEQ \esper[B_s B_t]=\half (|s|^{2\alpha}+|t|^{2\alpha}-|t-s|^{2\alpha}), \EEQ
is a natural generalization in the class of Gaussian processes of
the usual Brownian motion (which is the case $\alpha=\half$), in the sense that it exhibits two fundamental properties shared with Brownian motion, namely,
it has stationary increments, viz. $\esper[(B_t-B_s)(B_u-B_v)]=\esper[(B_{t+a}-B_{s+a})(B_{u+a}-B_{v+a})]$ for
every $a,s,t,u,v\in\R$, and it is self-similar, viz. 
\BEQ 
\forall \lambda>0, \quad (B_{\lambda t}, t\in\R) \overset{(law)}{=} (\lambda^{\alpha} B_t,
t\in\R).
\EEQ 
One may also define  a $d$-dimensional 
vector Gaussian process (called: {\it $d$-dimensional fractional Brownian motion}) by setting $B_t=(B_t(1),\ldots,B_t(d))$ where $(B_t(i),t\in\R)_{i=1,\ldots,d}$ are $d$ independent (scalar) fractional Brownian motions.

Its theoretical interest lies in particular in the fact that it is (up to normalization) the only Gaussian process satisfying these two properties.

A standard application of Kolmogorov's theorem shows that fBm has a version with
$\alpha^-$-H\"older continuous (i.e.  $\kappa$-H\"older continuous for every $\kappa<\alpha$) paths.
In particular, fBm with small Hurst parameter $\alpha$ is a natural, simple  model for continuous but
very irregular processes.

There has been a widespread interest during the past ten years in constructing a stochastic integration theory
with respect to fBm and solving stochastic differential equations driven by fBm, see for instance
\cite{LLQ02,GraNouRusVal04,CheNua05,RusVal93,RusVal00}. The multi-dimensional case is
very different from the one-dimensional case. When one tries to integrate for instance a stochastic differential
equation driven by a two-dimensional fBm $B=(B(1),B(2))$ by using any kind of Picard iteration scheme, one
encounters very soon the problem of defining the L\'evy area of $B$ which is the antisymmetric part
of ${\cal A}_{ts}:=\int_s^t dB_{t_1}(1) \int_s^{t_1} dB_{t_2}(2)$. This is the simplest occurrence
of iterated integrals $\vec{B}^k_{ts}(i_1,\ldots,i_k):=\int_s^t dB_{t_1}(i_1)\ldots \int_s^{t_{k-1}} dB_{t_k}(i_k)$,
$i_1,\ldots,i_k\le d$
for $d$-dimensional fBm $B=(B(1),\ldots,B(d))$ 
which lie at the heart of the rough path theory due to T. Lyons, see \cite{Lyo98,LyoQia02}. 
An alternative construction has been given by M. Gubinelli in \cite{Gu} under the name of 'algebraic
rough path theory', which we now propose to 
describe  briefly.

 Assume
$\Gamma_t=(\Gamma_t(1),\ldots,\Gamma_t(d))$ is some non-smooth $d$-dimensional path
which is $\alpha$-H\"older continuous.  Integrals such as
$\int f_1(\Gamma_t)d\Gamma_t(1)+\ldots+f_d(\Gamma_t)d\Gamma_t(d)$ do not make sense a priori because
$\Gamma$ is not differentiable (Young's integral \cite{Lej} works for $\alpha>\half$ but not beyond). In order to define
the integration of a differential form along $\Gamma$, it is enough to define a
 {\it geometric rough path} $(\vec{\Gamma}^{1},\ldots,\vec{\Gamma}^{\lfloor 1/\alpha \rfloor})$ lying above $\Gamma$,  
 $\lfloor  1/\alpha \rfloor$=entire part of $1/\alpha$,
  where $\vec{\Gamma}^{1}_{ts}=(\del\Gamma)_{ts}:=\Gamma_t-\Gamma_s$ is the {\it increment}
of $\Gamma$ between $s$ and $t$, and
each
$\vec{\Gamma}^k=(\vec{\Gamma}^k(i_1,\ldots,i_k))_{1\le i_1,\ldots,i_k\le d}$, $k\ge 2$ 
is a {\it substitute} for the iterated integrals $\int_s^t d\Gamma_{t_1}(i_1)\int_s^{t_1} d\Gamma_{t_2}(i_2)
\ldots \int_{s}^{t_{k-1}} d\Gamma_{t_k}(i_k)$ with the following three properties:

\begin{itemize}
\item[(i)] (\it H\"older continuity)
 each component of $\vec{\Gamma}^k$ is $k\alpha^-$-H\"older continuous, that is
to say, $k\kappa$-H\"older for every $\kappa<\alpha$;

\item[(ii)] ({\it multiplicativity}) letting $\del{\bf\Gamma}^k_{tus}:=
{\bf\Gamma}_{ts}^k-{\bf\Gamma}^k_{tu}-{\bf \Gamma}^k_{us}$, one requires
\BEQ
 \del\vec{\Gamma}^k_{tus}(i_1,\ldots,i_k) = \sum_{k_1+k_2=k} \vec{\Gamma}_{tu}^{k_1}(i_1,\ldots,i_{k_1}) \vec{\Gamma}_{us}^{k_2}(i_{k_1+1},\ldots,i_k); \label{eq:0:x}
\EEQ

\item[(iii)] ({\it geometricity}) 
\BEQ  {\bf \Gamma}^{n_1}_{ts}(i_1,\ldots,i_{n_1}) {\bf \Gamma}^{n_2}_{ts}(j_1,\ldots,j_{n_2}) 
 = 
\sum_{\vec{k}\in {\mathrm{Sh}}(\vec{i},\vec{j})} {\bf \Gamma}^{n_1+n_2}_{ts}(k_1,\ldots,k_{n_1+n_2})  \label{eq:0:geo} \EEQ
where ${\mathrm{Sh}}(\vec{i},\vec{j})$ is the subset of permutations of $i_1,\ldots,i_{n_1},j_1,\ldots,j_{n_2}$
which do not change the orderings of $(i_1,\ldots,i_{n_1})$ and $(j_1,\ldots,j_{n_2})$.
\end{itemize}

The multiplicativity property implies in particular the following identity for the
(non anti-symmetrized)  L\'evy area:
\BEQ {\cal A}_{ts}={\cal A}_{tu}+{\cal A}_{us}+ (B_t(1)-B_u(1))(B_u(2)-B_s(2)) \label{eq:0:mult} \EEQ
while the geometric property implies
\BEA && \int_s^t dB_{t_1}(1)\int_s^{t_1} dB_{t_2}(2)+\int_s^t dB_{t_2}(2)\int_s^{t_2} dB_{t_1}(1) \nonumber\\
&&=
\left(\int_s^t dB_{t_1}(1)\right) \left(\int_s^t dB_{t_2}(2)\right)=(B_t(1)-B_s(1))(B_t(2)-B_s(2)). \nonumber\\
\EEA

 Then there
is a standard procedure which allows to define out of these data iterated integrals of any order and
to solve differential equations driven by $\Gamma$.

The multiplicativity property (\ref{eq:0:x}) and the geometric property (\ref{eq:0:geo}) are satisfied by smooth paths, as can be checked by direct computation. So the most
natural way to construct such a multiplicative functional is to start from some smooth approximation
$\Gamma^{\eta}$, $\eta\overset{>}{\to} 0$ of $\Gamma$ such that each iterated
 integral $\vec{\Gamma}^{k,\eta}_{ts}(i_1,\ldots,i_k)$, $k\le \lfloor 1/\alpha \rfloor$ converges in the 
$k\kappa$-H\"older norm for every $\kappa<\alpha$.

This general scheme has been applied to fBm in a paper by L. Coutin and Z. Qian \cite{CQ02} and later in a paper
by the author \cite{Unt08}, using different schemes of approximation of $B$ by $B^{\eta}$ with
$\eta\to 0$.
 In both cases, the variance of the L\'evy area has been proved
to diverge in the limit $\eta\to 0$ when $\alpha\le 1/4$.

The approach developed in \cite{Unt08} makes use of a complex-analytic process $\Gamma$ defined
on the upper half-plane $\Pi^+=\{z=x+\II y\ | \ y>0\}$, called {\it $\Gamma$-process} or better
{\it analytic fractional Brownian motion} (afBm for short) \cite{TinUnt08}. Fractional Brownian motion $B_t$ appears
as the {\it real part} of the {\it boundary value} of $\Gamma_z$ when $\Im z\overset{>}{\to} 0$.
A natural approximation of $B_t$ is then obtained by considering 
\BEQ B_t^{\eta}:=\Gamma_{t+\II\eta}+\overline{\Gamma_{t+\II\eta}}=2\Re \Gamma_{t+\II\eta} \EEQ
 for $\eta\overset{>}{\to} 0.$ We show in subsection 3.1 that
$B^{\eta}$ may be written as a Fourier integral,
\BEQ B^{\eta}_t= c_{\alpha} \int_{\R} e^{-\eta|\xi|} |\xi|^{\half-\alpha} 
 \frac{e^{\II t\xi}-1}{\II \xi}
\ W(d\xi)  \label{eq:0:B-eta-Fourier} \EEQ
for some constant $c_{\alpha}$,
where $(W(\xi),\xi\ge 0)$ is a standard complex Brownian motion extended to $\R$
by setting $W(-\xi)=-\bar{W}(\xi)$, $\xi\ge 0$. When $\eta\to 0$, one retrieves
the well-known harmonizable representation of $B$ \cite{SamoTaq}.

The so-called {\it analytic iterated integrals} 
$$\int_s^t f_1(z_1) d\Gamma_{z_1}(1)\int_s^{z_1} f_2(z_2)d\Gamma_{z_2}(2)\ldots 
\int_s^{z_{d-1}} f_d(z_d) d\Gamma_{z_d}(d)$$
(where $f_1,\ldots,f_d$ are analytic functions), defined a priori for $s,t\in\Pi^+$ by integrating
over complex paths wholly contained in $\Pi^+$, converge to a finite limit when $\Im s,\Im t\to 0$
\cite{Unt08}, which is the starting point for the construction of a rough path associated to
$\Gamma$ \cite{TinUnt08}. The main tool for proving this kind of results is analytic continuation.

Computing iterated integrals
associated to $B_t=2\lim_{\eta\to 0} \Re\Gamma_{t+\II\eta}$
 instead of $\Gamma$ yields analytic iterated
 integrals, together with mixed integrals such as for instance
$\int_s^t d\Gamma_{z_1}(1)\int_s^{z_1} \overline{d\Gamma_{z_2}(2)}$. For these
 the analytic continuation method may no longer be applied because Cauchy's formula fails to hold,
and the above quantities may be shown to diverge when $\Re s,
\Re t\to 0$, see \cite{Unt08,Unt08b}.

\bigskip

\bigskip

Let us explain first how to define a L\'evy area for $B$. Proofs (as well as a sketch of
the Fourier normal ordering method for general iterated integrals) may be found in 
\cite{Unt09ter}.
 As mentioned before, the {\it uncorrected area}  ${\cal A}^{\eta}_{ts}:=\int_{s}^{t}
dB^{\eta}_{u_1}(1)\int_{s}^{u_1} dB^{\eta}_{u_2}(2)$ diverges when $\eta\to 0^+$.
The  idea is now to find some {\it increment counterterm} $(\del Z^{\eta})_{ts}=Z^{\eta}_{t}-Z^{\eta}_{s}$ such that
the {\it regularized area} ${\cal R}{\cal A}^{\eta}_{ts}:={\cal A}^{\eta}_{ts}-(\del Z^{\eta})_{ts}$ converges when $\eta\to 0^+$. Note that the multiplicativity  property (\ref{eq:0:mult}) holds for
${\cal R}{\cal A}^{\eta}$ as well as for ${\cal A}^{\eta}$ since $(\del Z^{\eta})_{ts}=
(\del Z^{\eta})_{t u}+(\del Z^{\eta})_{us}$. This counterterm $Z^{\eta}
$ may be found by using a suitable
decomposition of ${\cal A}^{\eta}_{ts}$ into the sum of :

-- an {\it increment term}, $(\del G^{\eta})_{ts}$;

-- a {\it boundary term} denoted by  ${\cal A}^{\eta}_{ts}(\partial)$.

The simplest idea one could think of would be to set 
\BEQ (\del G^{\eta})_{ts}=\int_{s}^{t}
dB^{\eta}_{u_1}(1)  B^{\eta}_{u_1}(2), \label{eq:0:G1}\EEQ
 and 
\BEQ {\cal A}^{\eta}_{ts}(\partial)=-\int_{s}^{t} dB^{\eta}_{u_1}(1)\ .\
B^{\eta}_{s}(2)=-B^{\eta}_{s}(2)(B^{\eta}_{t}(1)-B^{\eta}_{s}(1)).
\label{eq:0:del1} \EEQ
Alternatively, rewriting ${\cal A}^{\eta}_{ts}$ as $\int_{s}^{t} dB^{\eta}_{u_2}(2)
 \int_{u_2}^{t} dB^{\eta}_{u_1}(1)$, one may equivalently set

\BEQ (\del G^{\eta})_{ts}=-\int_{s}^{t}
dB^{\eta}_{u_2}(2) B^{\eta}_{u_2}(1) \label{eq:0:G2} \EEQ  and 
\BEQ {\cal A}^{\eta}_{ts}(\partial)=\int_{s}^{t} dB^{\eta}_{u_2}(2)\ .\
B^{\eta}_{t}(1)=B^{\eta}_{t}(1)(B^{\eta}_{t}(2)-B^{\eta}_{s}(2)).
\label{eq:0:del2} \EEQ

Now $\del G^{\eta}$ diverges when $\eta\to 0^+$,
 but since it is an increment, it may be discarded (i.e. it
might be used as a counterterm). The problem is, ${\cal A}^{\eta}_{ts}(\partial)$ converges
when $\eta\to 0^+$  in the 
$\kappa$-H\"older norm for every $\kappa<\alpha$, but not in the   $2\kappa$-H\"older norm (which is of course
well-known and may be seen as the starting point for
rough path theory). 

It turns out that a slight adaptation of this poor idea gives the solution. Decompose ${\cal A}^{\eta}_{ts}$ into a double integral in the Fourier coordinates
$\xi_1,\xi_2$ using (\ref{eq:0:B-eta-Fourier}). Use the
first increment/boundary decomposition (\ref{eq:0:G1},\ref{eq:0:del1}) 
for all indices $|\xi_1|\le |\xi_2|$, and the second one (\ref{eq:0:G2},\ref{eq:0:del2}) if 
$|\xi_1|>|\xi_2|$. Then
${\cal A}^{\eta}_{ts}(\partial)$, defined as the sum of  two contributions, one coming
from (\ref{eq:0:del1}) and the other from (\ref{eq:0:del2}),  {\em does converge} in the  
 $2\kappa$-H\"older norm when $\eta\to 0^+$, for every $\kappa<\alpha$.

As for the increment term $\del G^{\eta}$,
 defined similarly as the sum of two contributions coming from 
(\ref{eq:0:G1}) and (\ref{eq:0:G2}), it  diverges as soon as $\alpha\le 1/4$,
 but may be discarded at will. Actually we use in this
article a {\it minimal regularization scheme}: only the close-to-diagonal (i.e. $\xi_1/\xi_2\approx -1$)
 terms in the double integral
defining $\del G^{\eta}$ make it diverge. Summing over an appropriate subset, e.g.  $-\xi_1\not\in[\xi_2/2,2\xi_2]$ yields
an increment which converges (for {\it every} $\alpha\in(0,\half)$)
 when $\eta\to 0$ in the  $2\kappa$-H\"older norm for every $\kappa<\alpha$. 

Let $\alpha<1/4$. As noted in \cite{Unt08b}, the uncorrected L\'evy area ${\cal A}^{\eta}$ of the regularized process $B^{\eta}$ converges in law
to a  Brownian motion when $\eta\to 0^+$ after a rescaling by the factor $\eta^{\half(1-4\alpha)}$. In the latter 
article, the following question was raised: is it possible to define a counterterm $X^{\eta}$ living 
on the same probability
space as fBm, 
such that (i) the rescaled process $\eta^{\half(1-4\alpha)}X^{\eta}$ converges in law to Brownian motion;
 (ii) $(B^{\eta},{\cal A}^{\eta}- X^{\eta})$ is a multiplicative or almost multiplicative functional in the sense
of \cite{Lej}, Definition 7.1; (iii) ${\cal A}^{\eta}-X^{\eta}$ converges in the $2\kappa$-H\"older norm for
every $\kappa<\alpha$ when $\eta\to 0$ ? The 
  counterterm $X^{\eta}:={\cal A}^{\eta}-{\cal R} {\cal A}^{\eta}$ gives a solution to this problem.

\bigskip

The above ideas have a suitable generalization to  iterated integrals\\ 
$\int dB(i_1)\ldots \int dB(i_n)$ of order $n\ge 3$. There is one
more difficulty though: decomposing $(B^{\eta})'_{u_j}(i_j)$ into
 $c_{\alpha}\int dW_{\xi_j}(i_j) e^{\II u_j\xi_j} e^{-\eta|\xi_j|} |\xi_j|^{\half-\alpha}$,
 an extension
of the first increment/boundary decomposition (\ref{eq:0:G1}, \ref{eq:0:del1}), together with
a suitable regularization scheme,  yield the correct H\"older estimate {\em provided} $|\xi_1|\le \ldots\le |\xi_n|$. What should one
do then if $|\xi_{\sigma(1)}|\le \ldots\le |\xi_{\sigma(n)}|$ for some permutation $\sigma$
 instead ? The idea
is to permute the order of integration by using Fubini's theorem, and write $\int_{s}^{t}
 dB^{\eta}_{u_1}(i_1)\ldots \int_{s}^{u_{n-1}} dB^{\eta}_{u_n}(i_n)$ as
some {\em iterated tree integral} $\int dB^{\eta}_{u_1}(i_{\sigma(1)})\ldots \int dB^{\eta}_{u_n}
(i_{\sigma(n)})$.  The integration
domain, in the general case, becomes a little involved, and necessitates the introduction of combinatorial
tools on trees, such as admissible cuts for instance.  The underlying structures
are those of the  Hopf algebra of decorated rooted trees \cite{ConKre98,ConKre00}
  (as already noted in \cite{Kre99} or \cite{Gu2}), and of the Hopf shuffle algebra \cite{Mur1,Mur2}. The
proof of the multiplicative and of the geometric properties for the regularized rough
path, as well as the Hopf algebraic reinterpretation,  are to be found in \cite{Unt09bis}.
The general idea (see subsection 2.5 for more details) is that the fundamental objects
are {\em skeleton integrals} (a particular type of tree integrals) defined in subsection 2.1,
and that {\em any} regularization of the skeleton integrals (possibly even trivial)
 yielding finite quantities with the  correct H\"older regularity produces a regularized
rough path, which implies a large degree of arbitrariness in the definition.
 The idea  of cancelling singularities by building iteratively counterterms, originated from the Bogolioubov-Hepp-Parasiuk-Zimmermann (BPHZ) procedure
 for renormalizing Feynmann diagrams in quantum field theory \cite{Hepp}, mathematically
formalized in terms of Hopf algebras by A. Connes and D. Kreimer, has been applied during
the last decade in a variety of contexts ranging from numerical methods to quantum chromodynamics or 
multi-zeta functions, see for instance \cite{Kre99,Mur2,Wal00}. We plan to such a (less
arbitrary)   construction in the near future (see discussion at the end of subsection 2.5).

\bigskip

\bigskip
\bigskip

The main result of the paper may be stated as follows.

\begin{Theorem} \label{th:0}

{\it Let  $B=(B(1),\ldots,B(d))$  be a $d$-dimensional
fBm of Hurst index $\alpha\in(0,1)$, defined via the  harmonizable representation,
with the associated family of approximations $B^{\eta}$, $\eta>0$ living in the same 
probability space, see eq. (\ref{eq:0:B-eta-Fourier}). 
Then there exists a rough path
$({\cal R}{\bf B}^{1,\eta}=\del B^{\eta},\ldots,{\cal R}{\bf B}^{\lfloor 1/\alpha\rfloor,\eta})$ over
$B^{\eta}$ $(\eta>0)$, living in the chaos of order
$1,\ldots,\lfloor 1/\alpha\rfloor$ of $B$, 
satisfying properties (ii) (multiplicative property) and (iii) (geometric property) of
the Introduction, together with the following estimates:

\begin{itemize}
\item[(uniform H\"older estimate)] There exists a constant $C>0$ such that, for every $s,t\in\R$ and
$\eta>0$,
  $$\esper |{\cal R}{\bf B}^{n,\eta}_{ts}(i_1,\ldots,i_n)|^2 \le C|t-s|^{2n\alpha};$$
\item[(rate of convergence)] there exists a constant $C>0$ such that, for every $s,t\in\R$
and $\eta_1,\eta_2>0$,
   $$\esper |{\cal R}{\bf B}^{n,\eta_1}_{ts}(i_1,\ldots,i_n)-{\cal R}{\bf B}^{n,\eta_2}_{ts}(i_1,\ldots,i_n)|^2 \le C|\eta_1-\eta_2|^{2\alpha}.$$
\end{itemize}
}

\end{Theorem}

These results imply the existence of an explicit rough path ${\cal R}{\bf B}$ over $B$,
obtained as the limit of ${\cal R}{\bf B}^{\eta}$ when $\eta\to 0$.

\bigskip

\bigskip

Here is an outline of the article.
We first recall briefly some definitions and preliminary results on algebraic
rough path theory in Section 1, which show in particular that Theorem \ref{th:0}
implies the convergence of ${\cal R}{\bf B}^{\eta}$ to a rough path
${\cal R}{\bf B}$ over fractional Brownian motion $B$ when $\eta\to 0$.  Section 2
is dedicated to tree combinatorics and to the introduction of  quite general
regularization schemes for the iterated integrals of an arbitrary smooth path $\Gamma$. The
proof of the multiplicative and geometric properties are to be found in \cite{Unt09bis}
and are not reproduced here.
We apply a suitable regularization scheme to the construction of the regularized  rough path 
${\cal R}{\bf B}^{\eta}$  in section 3, and
prove the H\"older and rate of convergence estimates of Theorem \ref{th:0} for
the iterated integrals ${\cal R}{\bf B}^{n,\eta}(i_1,\ldots,i_n)$ with distinct
indices,
$i_1\not=\ldots\not=i_n$. We conclude in Section 4 by  showing how to extend these
results to coinciding indices, and introducing a new, real-valued, two-dimensional Gaussian  process
which we call {\it two-dimensional antisymmetric fractional Brownian motion}, to which the
above construction extends naturally.

\bigskip

{\bf Notations.} The group of permutations of $\{1,\ldots,n\}$ will be denoted by
$\Sigma_n$. The Fourier transform is ${\cal F}:f\to {\cal F}f(\xi)=\frac{1}{\sqrt{2\pi}}\int
f(x) e^{-\II x\xi} dx$.
If $|a|\le C|b|$ for some constant $C$ ($a$ and $b$ depending on some arbitrary set of parameters),
then we shall write $|a|\lesssim |b|$.


\section{The analysis of rough paths}


The present section will be very sketchy since the objects and results needed in this work have alread been
presented in great details in \cite{TinUnt08}. The fundational paper on the subject of algebraic rough path
theory is due to M. Gubinelli \cite{Gu}, see also \cite{Gu2} for more details in the case $\alpha<1/3$. Let us recall briefly the
original problem motivating the introduction of rough paths. Let $\Gamma:\R\to\R^d$ be some fixed
irregular (i.e. not differentiable) path, say $\kappa$-H\"older, and $f:\R\to\R^d$ some function which is also irregular (mainly because
one wants to consider functions $f$ obtained as a composition $g\circ \Gamma$ where $g:\R^d\to\R^d$ is regular).
 Can one define the integral $\int f_x d\Gamma_x$ ? The answer depends on the H\"older regularity of $f$ and $\Gamma$.
Assuming $f$ is $\gamma$-H\"older with $\kappa+\gamma>1$, then one may define the so-called {\it Young integral} \cite{Lej}
$\int_s^t f_x d\Gamma_x$ as the Riemann sum type limit $\lim_{|\Pi|\to 0} \sum_{\{t_j\}\in\Pi}
f_{t_i}(\Gamma_{t_{i+1}}-\Gamma_{t_i})$, where $\Pi=\{s=t_0<\ldots<t_n=t\}$ is a partition of $[s,t]$ with
mesh $|\Pi|$ going to $0$. Then the  resulting path $Y_t-Y_s:=\int_s^t f_x d\Gamma_x$ has the same regularity
as $\Gamma$. If $\kappa+\gamma\le 1$ instead, this is no more possible in general. One way out of this problem,
giving at the same time a coherent way to solve differential equations driven by $\Gamma$, is to define a class of 
{\it $\Gamma$-controlled paths} ${\cal Q}$, such that the above integration problem may be solved uniquely in
this class
by a formula generalizing the above Riemann sums, in which formal iterated integrals ${\bf \Gamma}^n(i_1,\ldots,i_n)$
of $\Gamma$ appear as in the Introduction.

\begin{Definition}[H\"older spaces]

{Let $\kappa\in(0,1)$ and $T>0$ fixed.
\begin{itemize}

\item[(i)] Let $C_1^{\kappa}=C_1^{\kappa}([0,T],\C)$ be the space of complex-valued $\kappa$-H\"older functions  $f$ in one variable
with (semi-)norm $||f||_{\kappa}=\sup_{s,t\in[0,T]} \frac{|f(t)-f(s)|}{|t-s|^{\kappa}}.$
\item[(ii)] Let $C_2^{\kappa}=C_2^{\kappa}([0,T],\C)$ be the space of complex-valued functions $f=f_{t_1,t_2}$
 of two variables vanishing on the diagonal $t_1=t_2$, such that $||f||_{2,\kappa}<\infty$, where
$|| \ .\  ||_{2,\kappa}$ is the following norm:
\BEQ ||f||_{2,\kappa}=\sup_{s,t\in[0,T]} \frac{|f_{t_1,t_2}|}{|t-s|^{\kappa}}.\EEQ
\item[(iii)]  Let $C_3^{\kappa}=C_3^{\kappa}([0,T],\C)$ be the space of complex-valued functions $f=f_{t_1,t_2,t_3}$
 of three variables vanishing on the subset $\{t_1=t_2\}\cup\{t_2=t_3\}\cup\{t_1=t_3\}$,
 such that $||f||_{3,\kappa}<\infty$ for some  generalized H\"older
semi-norm $|| \ .\ ||_{3,\kappa}$ defined for instance in \cite{Gu},
section 2.1.
\end{itemize}
}

\end{Definition}

\begin{Definition}[increments]

{\it 
\begin{itemize}
\item[(i)] Let $f$ be a function of one variable: then the increment of $f$, denoted by $\del f$, is
$(\del f)_{ts}:=f(t)-f(s).$
\item[(ii)] Let $f=f_{ts}$ be a function of two variables: then we define
\BEQ (\del f)_{tus}:= f_{ts}- f_{tu}-f_{us}.\EEQ
\end{itemize}
Note that $\del\circ\del(f)=0$ if $f$ is a function of one variable.
}
\end{Definition}

 Let $\Gamma=(\Gamma(1),\ldots,\Gamma(d)):[0,T]\to\R^d$ be a $\kappa$-H\"older path, and 
$({\bf\Gamma}^1_{ts}(i_1):=\Gamma_t(i_1)-\Gamma_s(i_1),{\bf\Gamma}^2_{ts}(i_1,i_2),\ldots,{\bf\Gamma}^{\lfloor
1/\kappa\rfloor}_{ts}(i_1,\ldots,i_{\lfloor 1/\kappa\rfloor}) )$ be a rough path lying above $\Gamma$, satisfying
properties (i) (H\"older property), (ii) (multiplicativity property)  and (iii) (geometricity property) of 
the Introduction. 

\begin{Definition}[controlled paths]

{\it Let $z=(z(1),\ldots,z(d))\in C_1^{\kappa}$ for some $\kappa<\alpha$ and $N=\lfloor 1/\kappa\rfloor+1$. Then $z$ is called a ($\Gamma$-)controlled path if its increments can be decomposed into
\BEQ \del z(i)=\sum_{n=1}^{N} \sum_{(i_1,\ldots,i_n)}  {\bf \Gamma}^n(i_1,\ldots,i_n).f^n(i_1,\ldots,i_n;i)+g^0(i)\EEQ
for some remainders $g^0(i)\in C_2^{N\kappa}$ and some paths $f^n(i_1,\ldots,i_n;i)\in (C_1^{\kappa})^n$ such that
\BEA &&  \del f^n(i_1,\ldots,i_n;i)= \nonumber\\
&&\quad \sum_{l=1}^{N-1-n}\sum_{(j_1,\ldots,j_l)} {\bf \Gamma}^l(j_1,\ldots,j_l). f^{l+n}(j_1,\ldots,j_l,i_1,\ldots,i_n;i)
+g^n(i_1,\ldots,i_n;i), \quad n=1,\ldots,N \nonumber\\ \EEA
for some remainder terms $g^n(i_1,\ldots,i_n;i)\in C_2^{(N-n)\kappa}$.

We denote by ${\cal Q}_{\kappa}$ the space of all such paths, and by ${\cal Q}_{\alpha^-}$ the intersection
$\cap_{\kappa<\alpha} {\cal Q}_{\kappa}.$
}

\end{Definition}

We may now state the main result.

\begin{Proposition}[see \cite{Gu2}, Theorem 8.5, or \cite{TinUnt08}, Proposition 3.1]

{\it Let $z\in{\cal Q}_{\alpha^-}$. Then the limit
\BEQ \int_s^t z_x d\Gamma_x:=\lim_{|\Pi|\to 0} \sum_{k=0}^n \sum_{i=1}^d  \left[ \del X_{t_{k+1},t_k}(i) z_{t_k}(i)
+\sum_{n=1}^{N-1}\sum_{(i_1,\ldots,i_{n})} {\bf \Gamma}^{n+1}_{t_{k+1},t_k}(i_1,\ldots,i_n,i)\zeta^n_{t_k}(i_1,\ldots,i_n;i)
 \right] \EEQ
exists in the space ${\cal Q}_{\alpha^-}$.
}

\end{Proposition}

Assume $\Gamma$ is a centered Gaussian process, and $\Gamma^{\eta}$ a family of Gaussian
approximations of $\Gamma$ living in its first chaos. Then
the Proposition below gives very convenient moment conditions for a family of
 rough paths $(\Gamma^{\eta},{\bf\Gamma}^{2,\eta},\ldots,{\bf\Gamma}^{\lfloor 1/\kappa\rfloor,\eta})$
to converge in the right H\"older norms when $\eta\to 0$, thereby
   defining a rough path  above  $\Gamma$.

\begin{Proposition}

{\it Let $\Gamma$ be a $d$-dimensional centered Gaussian process admitting a version with a.s. $\alpha^-$-H\"older paths. Let
$N=\lfloor 1/\alpha\rfloor.$ 
Assume:
\begin{enumerate}
\item  there exists a family $\Gamma^{\eta}$, $\eta\to 0^+$ of Gaussian processes living in 
the first chaos  of $\Gamma$ and an overall constant $C$ such that 
\begin{itemize}
\item[(i)] \BEQ \esper |\Gamma^{\eta}_t-\Gamma_s^{\eta}|^2 \le C|t-s|^{2\alpha}; \label{eq:1:i} \EEQ
\item[(ii)] \BEQ \esper |\Gamma^{\eta}_t-\Gamma^{\eps}_t|^2 \le C |\eps-\eta|^{2\alpha}; \label{eq:1:ii} \EEQ
\item[(iii)] $\forall t\in[0,T]$, $\Gamma_t^{\eta}\overset{L^2}{\to} \Gamma_t$ when $\eta\to 0$;
\end{itemize}
\item there exists a truncated multiplicative functional $({\bf \Gamma}_{ts}^{1,\eta}=\Gamma_t^{\eta}-\Gamma_s^{\eta},
{\bf\Gamma}_{ts}^{2,\eta},\ldots,{\bf\Gamma}_{ts}^{N,\eta})$ lying above $\Gamma^{\eta}$ and living
in the $n$-th chaos of $\Gamma$, $n=1\ldots,N$,  such that,
for every $2\le k\le N$,
\begin{itemize}
\item[(i)] \BEQ \esper |{\bf\Gamma}_{ts}^{k,\eta}|^2 \le C|t-s|^{2k\alpha}; \label{eq:1:ibis}\EEQ
\item[(ii)] \BEQ \esper |{\bf \Gamma}_{ts}^{k,\eps}-{\bf\Gamma}_{ts}^{k,\eta}|^2 \le C|\eps-\eta|^{2\alpha}.
\label{eq:1:iibis} \EEQ
\end{itemize}
\end{enumerate}

Then $({\bf\Gamma}^{1,\eta},\ldots,{\bf\Gamma}^{N,\eta})$ converges in
$L^2(\Omega;C_2^{\kappa}([0,T],\R^d)\times C_2^{2\kappa}([0,T],\R^{d^2}) \times\ldots\times C_2^{N\kappa}([0,T],\R^{d^N}))$  for every $\kappa<\alpha$ to a rough path $({\bf\Gamma}^1,\ldots,{\bf\Gamma}^N)$ lying above $\Gamma$.
}
\end{Proposition}

{\bf Short proof} (see \cite{TinUnt08}, Lemma 5.1, Lemma 5.2 and Prop. 5.4).
 The main ingredient is the Garsia-Rodemich-Rumsey (GRR for
short) lemma \cite{Ga} which
states that, if $f\in C_2^{\kappa}([0,T],\C)$,
\BEQ  ||f||_{2,\kappa} \le  C \left(
 ||\del f||_{3,\kappa}+ \left( \int_0^T \int_0^T \frac{|f_{vw}|^{2p}}{|w-v|^{2\kappa p+2}} \ dv\ dw\right)^{1/2p}
\right)  \EEQ
 for every $p\ge 1$.

Then properties (\ref{eq:1:i},\ref{eq:1:ibis})  imply by using the  GRR lemma for $p$ large enough,
 Jensen's inequality and the
equivalence of $L^p$-norms for processes living in a fixed Gaussian chaos 
\BEQ \esper ||{\bf\Gamma}^{k,\eta}||_{2,k\kappa} \lesssim \esper ||\del{\bf\Gamma}^{k,\eta}||_{3,k\kappa}+C.\EEQ
By using the multiplicative property (ii) in the Introduction and induction on $k$, $\esper ||\del{\bf\Gamma}^{k,\eta}||_{3,k\kappa}$ may in the same way be proved to be bounded by a constant.

On the other hand, properties (\ref{eq:1:i},\ref{eq:1:ii},\ref{eq:1:ibis},\ref{eq:1:iibis}), together
with the equivalence of $L^p$-norms,  imply (for every $\kappa<\alpha$)
\BEQ \esper |{\bf\Gamma}^{k,\eps}_{ts}-{\bf\Gamma}^{k,\eta}_{ts}|^2 \lesssim |t-s|^{2k\kappa} |\eps-\eta|^{2(\alpha-\kappa)}\EEQ
hence, by the same arguments,
\BEQ \esper ||{\bf\Gamma}^{k,\eps}-{\bf\Gamma}^{k,\eta}||_{2,k\kappa} \lesssim |\eps-\eta|^{\alpha-\kappa} \EEQ
which shows that ${\bf\Gamma}^{k,\eps}$ is a Cauchy sequence in $C_2^{k\kappa}([0,T],\R^{d^k})$. \hfill \eop


\section{Tree combinatorics and the Fourier normal ordering method}



\subsection{From iterated integrals to trees}


It was noted already long time ago \cite{But72} that iterated integrals could be encoded by trees. This
remark has been exploited in connection with the construction of the rough path
solution of (partial, stochastic) differential equations  in \cite{Gu2}.
The correspondence between trees and itegrated integrals goes simply as follows.

\begin{Definition}

A decorated rooted tree (to be drawn growing {\em up}) is a finite tree with a distinguished
vertex called {\em root} and edges oriented {\em downwards} (i.e. directed towards the root), such
that every vertex wears an integer label.

If $\T$ is a decorated rooted tree, we let $V(\T)$ be the set of its vertices (including the root),
and $\ell:V(\T)\to\N$ be its vertex labeling.

More generally, a decorated rooted {\em forest} is a finite set of decorated rooted trees. If $\T=\{\T_1,
\ldots,\T_l\}$ is a forest, then we shall write $\T$ as the formal commutative product $\T_1\ldots\T_l$.

\end{Definition}

\begin{Definition}  \label{def:2:connect}

Let $\T$ be a decorated rooted tree.

\begin{itemize}
\item Letting $v,w\in V(\T)$, we say that $v$ {\em connects directly to} $w$, and write
$v\to w$ or equivalently $w=v^-$, if $(v,w)$ is an edge oriented downwards from $v$
to $w$. (Note that $v^-$ exists and is unique except if $v$ is the root).
\item If $v_m\to v_{m-1}\to\ldots\to v_1$, then we shall write $v_m\twoheadrightarrow v_1$, and say that
$v_m$ {\em connects to} $v_1$. By definition, all vertices (except the root) connect to the root.
\item Let $(v_1,\ldots,v_{|V(\T)|})$ be an ordering of $V(\T)$. Assume that
$\left( v_i\twoheadrightarrow v_j\right)\Rightarrow\left(i>j\right)$ (in particular, $v_1$ is the root).
Then we shall say that the ordering is {\em compatible} with the {\em tree partial ordering} defined
by $\twoheadrightarrow$.
\end{itemize}

\end{Definition}

\begin{figure}[h]
  \centering
   \includegraphics[scale=0.35]{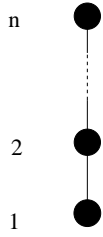}
   \caption{\small{Trunk tree.}}
  \label{Fig1}
\end{figure}

\begin{Definition} \label{def:2:it-int}

\begin{itemize}
\item[(i)]
Let  $\Gamma=(\Gamma(1),\ldots,\Gamma(d))$ be a $d$-dimensional smooth path, and $\T$ a decorated rooted
tree such that $\ell:V(\T)\to\{1,\ldots,d\}$. Then $I_{\T}(\Gamma):\R^2\to\R$ is the iterated integral
defined as
\BEQ [I_{\T}(\Gamma)]_{ts}:=\int_s^t d\Gamma_{x_{v_1}}(\ell(v_1))\int_s^{x_{v_2^-}} 
d\Gamma_{x_{v_2}}(\ell(v_2))
\ldots \int_s^{x_{v^-_{|V(\T)|}}} d\Gamma_{x_{v_{|V(\T)|}}}(\ell(v_{|V(\T)|}))
\EEQ
where $(v_1,\ldots,v_{|V(\T)|})$ is any ordering of $V(\T)$ compatible with the tree partial ordering.

In particular, if $\T$ is a trunk tree with $n$ vertices (see Fig. \ref{Fig1}) -- so that
the tree ordering is total -- we shall write
\BEQ I_{\T}(\Gamma)=I_n^{\ell}(\Gamma),\EEQ
where
\BEQ [I_n^{\ell}(\Gamma)]_{ts}:=\int_s^t d\Gamma_{x_1}(\ell(1)) \int_s^{x_1} d\Gamma_{x_2}(\ell(2))
\ldots \int_s^{x_{n-1}} d\Gamma_{x_{n}}(\ell(n)).
\EEQ

\item[(ii)] (generalization)
Assume $\T$ is a subtree of $\tilde{\T}$. Let $\mu$ be a Borel measure on $\R^{\tilde{\T}}$.
Then
\BEQ [I_{\tilde{\T}}(\mu)]_{ts}:=\int_s^t \int _s^{x_{v_1^-}}\ldots\int_s^{x_{v^-_{|V(\T)|}}}
\mu(dx_{v_1},\ldots,dx_{v_{|V(\T)|}}) \EEQ
is a measure on $\R^{\tilde{\T}\setminus\T}$.

\end{itemize}

\end{Definition}

Assume $\T=\tilde{\T}$ so $[I_{\tilde{\T}}(\mu)]_{ts}$ is a {\it number}. Then case (i)
may be seen as a particular case of case (ii) with $\mu=d\Gamma(\ell(v_1))\otimes\ldots
\otimes  d\Gamma(\ell(v_{|V(\T)|}))$. Conversely, case (ii) may be seen as a multilinear extension
of case (i), and will turn out to be useful later on for the regularization procedure. Note
however that (i) uses the labels of $\T$ while (ii) {\it doesn't}.

The above correspondence extends by (multi)linearity to the {\em algebra of decorated rooted 
trees} which we shall
 now introduce.

\begin{Definition}[algebra of decorated rooted trees]

\begin{itemize}
\item[(i)] Let $\cal T$ be the free commutative algebra over $\R$ generated by decorated rooted trees.
If $\T_1,\T_2,\ldots \T_l$ are (decorated rooted) trees, then the product  $\T_1\ldots\T_l$ is the
forest with connected components $\T_1,\ldots,\T_l$. 
\item[(ii)] Let $\T'=\sum_{l=1}^L m_l \T_l\in{\cal T}$, where $m_l\in\R$ and each $\T_l=\T_{l,1}\ldots
\T_{l,L(l)}$ is a forest with labels in the set $\{1,\ldots,d\}$, and $\Gamma$ be a smooth $d$-dimensional
path as above. Then
\BEQ [I_{\T'}(\Gamma)]_{ts}:=\sum_{l=1}^L m_l [I_{\T_{l,1}}(\Gamma)]_{ts}
 \ldots [I_{\T_{l,L(l)}}(\Gamma)]_{ts}. \EEQ
\end{itemize}

\end{Definition}

\medskip

Let us now rewrite these iterated integrals by using Fourier transform. 

\begin{Definition}[formal integral]

Let $f:\R\to\R$ be a smooth, compactly supported function such that ${\cal F}f(0)=0$.
Then the formal integral $\int^t f=-\int_t f$ of $f$ is defined as $\frac{1}{\sqrt{2\pi}}
\int_{-\infty}^{+\infty} ({\cal F}f)(\xi) \frac{e^{\II t\xi}}{\II \xi}\ d\xi.$

\end{Definition}

Formally one may write: 
\BEQ \int^t e^{\II x\xi} dx=\int_{\pm\II\infty}^t e^{\II x\xi} dx=\frac{e^{\II t\xi}}{\II\xi} \EEQ
(depending on the sign of $\xi$). The condition ${\cal F}f(0)=0$ prevents possible infra-red divergence
when $\xi\to 0$.

\medskip

The {\it skeleton integrals} defined below must be understood in a {\em formal} sense because of
the possible infra-red divergences.

\begin{Definition}[skeleton integrals]

\begin{itemize}
\item[(i)]

Let $\T$ be a tree with $\ell:\T\to\{1,\ldots,d\}$ and $\Gamma$ be a $d$-dimensional compactly
supported, smooth path. Let $(v_1,\ldots,
v_{|V(\T)|})$ be any ordering of $V(\T)$ compatible with the tree partial ordering. Then the
{\em skeleton integral} of $\Gamma$ along $\T$ is by definition
\BEQ [\Sk I_{\T}(\Gamma)]_t=\int^t d\Gamma_{x_{v_1}}(\ell(v_1))
\int^{x_{v_2^-}} d\Gamma_{x_2}(\ell(v_2))
\ldots \int^{x_{v^-_{|V(\T)|}}} d\Gamma_{x_{v_{|V(\T)|}}}(\ell(v_{|V(\T)|})).
\EEQ

\item[(ii)] (multilinear extension, see Definition \ref{def:2:it-int}) Assume $\T$ is
a subtree of $\tilde{\T}$, and $\mu$ a compactly supported Borel measure on $\R^{\tilde{\T}}$. Then
\BEQ [\Sk I_{\T}(\mu)]_t=\int^t \int^{x_{v_2^-}}\ldots\int^{x_{v^-_{|V(\T)|}}} \mu(dx_{v_1},\ldots,
dx_{v_{|V(\T)|}}) \EEQ
is a measure on $\R^{\tilde{\T}\setminus\T}$.

\end{itemize}

\end{Definition}

Formally again, $[\Sk I_{\T}(\Gamma)]_t$ may be seen as $[I_{\T}(\Gamma)]_{t,\pm\II\infty}$.
Note that (denoting by $\hat{\mu}$ the partial Fourier transform of $\mu$ with respect to $(x_v)_{v\in V(\T)}$),
the following equation holds,
\BEQ [\SkI_{\T}(\mu)]_t=(2\pi)^{-|V(\T)|/2} \langle \hat{\mu}, \left[
\SkI_{\T}\left( (x_v)_{v\in V(\T)} \to  e^{\II \sum_{v\in V(\T)} x_v \xi_v} \right) \right]_t \rangle.
\label{eq:2:reg-sk-it0}  \EEQ

\begin{Lemma} \label{lem:2:SkI}

The following formula holds:
\BEQ [\Sk I_{\T}(\Gamma)]_t=(\II\sqrt{2\pi})^{-|V(\T)|} \int\ldots\int_{\R^{\T}} \prod_{v\in V(\T)} d\xi_v
\ .\ e^{\II t\sum_{v\in V(\T)} \xi_v} \frac{\prod_{v\in V(\T)} {\cal F}(\Gamma'(\ell(v)))(\xi_v)}{\prod_{v\in V(\T)} (\xi_v+\sum_{w\twoheadrightarrow v} \xi_w)}.\EEQ

\end{Lemma}

{\bf Proof.} We use induction on $|V(\T)|$. After stripping the root of $\T$ (denoted by $0$) there remains
a forest $\T'=\T'_1\ldots\T'_J$, whose roots are the vertices directly connected to $0$. Assume
\BEQ [\Sk I_{\T'_j}(\Gamma)]_{x_0}=\int\ldots\int \prod_{v\in V(\T'_j)} d\xi_v \ .\ e^{\II x_0\sum_{v\in V(\T'_j)} \xi_v} F_j((\xi_v)_{v\in \T'_j}).\EEQ
Note that
\BEQ {\cal F} \left( \prod_{j=1}^J \Sk I_{\T'_j}(\Gamma) \right)(\xi)=\int_{\sum_{v\in V(\T)\setminus
\{0\}} \xi_v=\xi} \prod_{v\in V(\T)\setminus\{0\}} d\xi_v \prod_{j=1}^J F_j((\xi_v)_{v\in V(\T'_j)}).\EEQ

Then
\BEA
[\Sk I_{\T}(\Gamma)]_t&=& \int^t d\Gamma_{x_0}(\ell(0)) \prod_{j=1}^J [\Sk I_{\T'_j}(\Gamma)]_{x_0} 
\nonumber\\
&=& \frac{1}{\sqrt{2\pi}}\int_{-\infty}^{+\infty} \frac{d\xi}{\II\xi}  e^{\II t\xi} {\cal F}\left( \Gamma'(\ell(0))
\prod_{j=1}^J \Sk I_{\T'_j}(\Gamma) \right)(\xi) \nonumber\\
&=&  \frac{1}{\sqrt{2\pi}} \int_{-\infty}^{+\infty} d\xi_0 {\cal F}(\Gamma'(\ell(0)))(\xi_0) e^{\II t\xi_0} \ .\ \nonumber\\
&& \qquad \qquad \int_{-\infty}^{+\infty} \frac{d\xi}{\II\xi} e^{\II t(\xi-\xi_0)}
\int_{\sum_{v\in V(\T)\setminus\{0\}} \xi_v=\xi-\xi_0}  d\xi_v \prod_{j=1}^J F_j((\xi_v)_{v\in V(\T'_j)})
\nonumber\\ 
\EEA
hence the result. \hfill \eop

\bigskip

Skeleton integrals are the fundamental objects from which regularized rough paths will be constructed
in the next subsections.


\subsection{Coproduct structure and increment-boundary decomposition}


Consider for an example the trunk tree $\T^{\Id_n}$ (see subsection 2.4 for an explanation of the
notation)  with vertices $n\to n-1\to\ldots\to 1$ and
labels $\ell:\{1,\ldots,n\}\to\{1,\ldots,d\}$, and the associated iterated integral (assuming
$\Gamma=(\Gamma(1),\ldots,\Gamma(d))$ is a smooth path)
\BEQ [I_n^{\ell}(\Gamma)]_{ts}=
 [I_{\T^{\Id_n}}(\Gamma)]_{ts}=\int_s^t d\Gamma_{x_1}(\ell(1))\ldots\int_s^{x_{n-1}}d\Gamma_{x_{n}}
(\ell(n)).\EEQ

Cutting $\T^{\Id_n}$ at some vertex $v\in\{2,\ldots,n\}$ produces two trees, $L_v \T^{\Id_n}$ 
({\em left} or rather {\em bottom}
 part of $\T^{\Id_n}$) and $R_v\T^{\Id_n}$ ({\em right} or {\em top} part), with
respective vertex subsets $\{1,\ldots,v-1\}$ and $\{v,\ldots,n\}$. One should actually see the couple
$(L_v\T^{\Id_n},R_v\T^{\Id_n})$ as $L_v\T^{\Id_n}\otimes R_v\T^{\Id_n}$ sitting in the tensor
product algebra ${\cal T}\otimes{\cal T}$. Then  multiplicative property (ii) in the Introduction
reads
\BEQ [\del I_{\T^{\Id_n}}(\Gamma)]_{tus}=\sum_{v\in V(\T^{\Id_n})\setminus\{1\}} [I_{L_v\T^{\Id_n}}(\Gamma)]_{tu}
[I_{R_v \T^{\Id_n}}(\Gamma)]_{us}. \label{eq:treex0} \EEQ

On the other hand, one may rewrite $[I_{\T^{\Id_n}}(\Gamma)]_{ts}$ 
as the sum of the {\em increment term}
\BEA  [\del G]_{ts} & =\int^t d\Gamma_{x_1}(\ell(1))\int^{x_1} d\Gamma_{x_2}(\ell(2))
\ldots\int^{x_{n-1}}d\Gamma_{x_{n}}(\ell(n))\nonumber\\
& \qquad   - \int^s d\Gamma_{x_1}(\ell(1))\int^{x_1} d\Gamma_{x_2}(\ell(2))
\ldots\int^{x_{n-1}}d\Gamma_{x_{n}}(\ell(n)) \label{eq:2:incr-term} \nonumber\\  \EEA
and of the {\em boundary term}
\BEA &&  [I_{\T^{\Id_n}}(\Gamma)(\partial)]_{ts}=  -\sum_{n_1+n_2=n} \int_s^t d\Gamma_{x_1}(\ell(1))
\ldots \int_s^{x_{n_1-1}}d\Gamma_{x_{n_1}}(\ell(n_1)) \ .\ \nonumber\\
 && \qquad  .\ \int^s d\Gamma_{x_{n_1+1}}(\ell(n_1+1))
\int^{x_{n_1+1}}d\Gamma_{x_{n_1+2}}(\ell(n_1+2)) \ldots \int^{x_{n-1}} d\Gamma_{x_{n}}(\ell(n)). 
\label{eq:2:bdry-term} \nonumber\\
\EEA

The above decomposition is fairly obvious for $n=2$ (see Introduction) and obtained by easy induction for
general $n$. Thus  (using tree notation this time)
\BEQ [I_{\T^{\Id_n}}(\Gamma)]_{ts}=[\del \Sk I_{\T^{\Id_n}}]_{ts}-\sum_{v\in V(\T^{\Id_n})\setminus\{1\}}
[I_{L_v\T^{\Id_n}}(\Gamma)]_{ts} \ .\ [\Sk I_{R_v\T^{\Id_n}}(\Gamma)]_s. \label{eq:2:*} \EEQ

\medskip

The above considerations extend to arbitrary trees (or also forests) as follows.

\begin{Definition}[admissible cuts]

{\it
\begin{enumerate}
\item
 Let $\T$ be a tree, with set of vertices $V(\T)$ and root
denoted by  $0$.
If $\vec{v}=(v_1,\ldots,v_J)$, $J\ge 1$  is any totally disconnected subset of $V(\T)\setminus\{0\}$,
 i.e. $v_i\not\twoheadrightarrow
v_j$ for all $i,j=1,\ldots,J$, then we shall say that $\vec{v}$ is an {\em admissible cut} of $\T$, and
 write $\vec{v}\models V(\T)$. We let $R_{\vec{v}}\T$ be the sub-forest (or sub-tree if $J=1$) obtained by keeping
only the vertices above $\vec{v}$, i.e. $V(R_{\vec{v}}\T)=\vec{v}\cup\{w\in V(\T):\ \exists j=1,\ldots,J,
w\twoheadrightarrow v_j\}$, and $L_{\vec{v}}\T$ be the sub-tree obtained by keeping all other vertices.

\item Let $\T=\T_1\ldots\T_l$ be a forest, together with its decomposition into trees. Then an admissible
cut of $\T$ is a disjoint union $\vec{v}_1\cup\ldots\cup\vec{v}_l$, $\vec{v}_i\subset\T_i$, where $\vec{v}_i$
is either $\emptyset$, $\{0_i\}$ (root of $\T_i$) or an admissible cut of $\T_i$. By
definition, we let $L_{\vec{v}}\T=L_{\vec{v}_1}\T_1 \ldots L_{\vec{v}_l}\T_l$, $R_{\vec{v}}\T=
R_{\vec{v}_1}\T_1\ldots R_{\vec{v}_l}\T_l$ (if $\vec{v}_i=\emptyset$, resp. $\{0_i\}$, then
$(L_{\vec{v}_i}\T_i,R_{\vec{v}_i}\T_i):=(\T_i,\emptyset)$, resp. $(\emptyset,\T_i)$). 

We exclude by
convention the two trivial cuts $\emptyset\cup\ldots\cup\emptyset $ and $\{0_1\}\cup\ldots\cup\{0_l\}$.
\end{enumerate}
}
\label{def:6:admissible-cut}
\end{Definition}

\begin{figure}[h]
  \centering
   \includegraphics[scale=0.35]{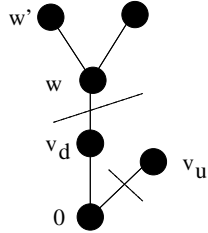}
   \caption{\small{Admissible cut.}}
  \label{Fig2}
\end{figure}

\begin{figure}[h]
  \centering
   \includegraphics[scale=0.35]{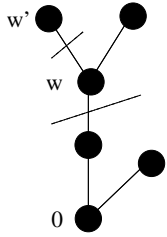}
   \caption{\small{Non-admissible cut.}}
  \label{Fig3}
\end{figure}

See Fig. \ref{Fig2} and \ref{Fig3}. Defining the co-product operation $\Del:
{\cal T}\to {\cal T}\otimes
{\cal T}$, $\T\to e \otimes \T+\T\otimes e+\sum_{\vec{v}\models V(\T)} L_{\vec{v}}\T\otimes
R_{\vec{v}}\T$ (where $e$ stands for the empty tree, which is the unit of the algebra)
yields a coalgebra structure on ${\cal T}$ which makes it (once the antipode -- which we do
not need here --  is defined) a Hopf algebra (see articles by A. Connes
and D. Kreimer \cite{ConKre98,ConKre00,ConKre01}). The convention is usuall to write $\vec{v}=c$ (cut), $L_{\vec{v}}\T=R^c (\T)$ (root part), $R_{\vec{v}}\T=P^c (\T)$ and $\Del(\T)=
e\otimes\T+\T\otimes e+\sum_c P^c(\T)\otimes R^c(\T)$ (note the inversion of the order
of the factors in the tensor product).

Eq. (\ref{eq:treex0}) extends to the general formula (called: {\it tree multiplicative property}), which
one can find in \cite{Kre99} or \cite{Gu2},
\BEQ [\del I_{\T}(\Gamma)]_{tus}=\sum_{\vec{v}\models V(\T)} [I_{L_{\vec{v}}\T}(\Gamma)]_{tu}
[I_{R_{\vec{v}} \T}(\Gamma)]_{us}, \label{eq:treex} \EEQ
satisfied by any regular path $\Gamma$ for any tree $\T$.

Letting formally $s=\pm\II\infty$ in eq. (\ref{eq:treex})  yields

\BEQ [I_{\T}(\Gamma)]_{tu}=[\del \Sk I_{\T}]_{tu}-\sum_{v\in V(\T)\setminus\{0\}}
[I_{L_v\T}(\Gamma)]_{tu} \ .\ [\Sk I_{R_v\T}(\Gamma)]_u. \label{eq:2:SkI} \EEQ

which generalizes eq. (\ref{eq:2:*}). Conversely, eq. (\ref{eq:2:SkI}) implies the tree multiplicative property eq. (\ref{eq:treex}), as shown in Lemma \ref{lem:2:reg} below.


\subsection{Regularization procedure}


\begin{Definition}[regularization procedure for skeleton integrals] \label{def:2:reg}

Let $\tilde{\T}=\{v_1<\ldots<v_{|\tilde{T}|}\}$ be a tree, $\T\subset\tilde{\T}$ a subtree, $\mu$ a compactly supported Borel measure on
$\R^{\tilde{\T}}$ such that $\supp\hat{\mu}\subset\{(\xi_1,\ldots,\xi_{|V(\tilde{T})|}) \ | \
|\xi_1|\le \ldots\le |\xi_{|V(\tilde{\T})|}\}$, and $D_{reg}\subset\R^{\T}$ a Borel subset.

The (formal) $D_{reg}$-{\em regularized skeleton integral} ${\cal R}\SkI_{\T}$
 is the linear mapping (see eq. (\ref{eq:2:reg-sk-it0}))
\BEQ \mu\to [{\cal R}\SkI_{\T}(\mu)]_{t}=(2\pi)^{-|V(\T)|/2} \langle \hat{\mu}, {\bf 1}_{D_{reg}}(\xi) \ .\
\left[ \SkI_{\T}\left( (x_v)_{v\in V(\T)} \to e^{\II \sum_{v\in V(\T)} x_v \xi_v} \right) \right]_t \rangle 
\label{eq:2:reg-sk-it} \EEQ
where $\hat{\mu}$ is  the partial Fourier transform of $\mu$ with respect to $(x_v)_{v\in V(\T)}$.

By assumption we shall only allow $D_{reg}=\R$ if $\T$ is a tree reduced to one vertex.

\end{Definition}

\begin{Lemma}[regularization] \label{lem:2:reg}
{\it

Let $\T=\T_1\ldots\T_l$ be a  forest, together with its tree decomposition.
Define by induction on $|V(\T)|$ the regularized integration operator  $\left[{\cal R}I_{\T}\right]_{ts}$ by
\BEQ \prod_{j=1}^l \left\{ 
\left[ \del {\cal R}\SkI_{ \T_j}\right]_{ts}-
\sum_{\vec{v}\models V(\T_j)}
\left[ {\cal R} I_{L_{\vec{v}}\T_j}\right]_{ts} \left[ {\cal R}\SkI_{ R_{\vec{v}}\T_j}\right]_{s}
\right\}    \label{eq:6:reg-int} \EEQ

Then $\left[{\cal R}I_{\T}\right]_{ts}$ satisfies the following tree multiplicative property:

\BEQ \left[ \del{\cal R}I_{\T}\right]_{t us}
=\sum_{\vec{v}\models V(\T)} \left[ {\cal R} I_{L_{\vec{v}} \T }
 \right]_{t u} \ .\ \left[ {\cal R} I_{R_{\vec{v}} \T} \right]_{us}. \label{eq:6:x}  \EEQ

By analogy with eq. (\ref{eq:2:incr-term}, \ref{eq:2:bdry-term}, \ref{eq:2:*}),
 $\left[ \del {\cal R}\SkI_{\T_j}\right]_{ts}$, resp.
$[{\cal R} I_{\T_j}(\partial)]_{ts}:= -
\sum_{\vec{v}\models V(\T_j)}
\left[ {\cal R} I_{L_{\vec{v}}\T_j}\right]_{ts} \left[ {\cal R}\SkI_{ R_{\vec{v}}\T_j}\right]_{s} $
 may be called the {\em increment},
resp. {\em boundary operators} associated to the tree $\T_j$.
}

\end{Lemma}

{\bf Remark.} By Definition \ref{def:2:reg}, the condition $[{\cal R}I_{\T}]_{ts}=[I_{\T}]_{ts}$
holds for a tree reduced to one vertex. This implies in the end that one has constructed a 
rough path over the {\em original} path $\Gamma$.

{\bf Proof.}  If the multiplicative property (\ref{eq:6:x}) holds for trees, then it holds automatically
for forests since $[{\cal R} I_{\T_1\ldots\T_l}]_{ts}$ is the product
$\prod_{j=1}^l [{\cal R}I_{\T_j}]_{ts}$. Hence we may assume that   $\T$ is a tree, say, with
$n$ vertices.
 Suppose (by induction) that the above multiplicative property
(\ref{eq:6:x}) holds for all  trees with $\le n-1$ vertices. Then
\BEA  && \left[ \del{\cal R}I_{\T}\right]_{t u s}  = 
    \sum_{\vec{v}\models V(\T)}  \left( -
\left[ \del {\cal R} I_{L_{\vec{v}}\T}\right]_{t u s}
                 \left[ {\cal R}\SkI_{R_{\vec{v}}\T}\right]_{s}+ 
 \left[ {\cal R} I_{L_{\vec{v}}\T}\right]_{t u}
                 \left[\del {\cal R}\SkI_{R_{\vec{v}}\T}\right]_{us} \right)  \nonumber\\
&& = \sum_{\vec{v}\models V(\T)} \sum_{\vec{w}\models V(L_{\vec{v}} \T)}  \left( - \left[ {\cal R}
I_{ L_{\vec{w}}\circ L_{\vec{v}}(\T)} \right]_{t u}
\left[ {\cal R} I_{R_{\vec{w}}\circ L_{\vec{v}}(\T)} \right]_{us} \left[{\cal R} \SkI_{R_{\vec{v}}\T}\right]_{s}
\right. \nonumber\\
&& \left. \qquad \qquad \qquad  
 +  \left[ {\cal R} I_{L_{\vec{v}}\T}\right]_{t u}
                 \left[ \del {\cal R}\SkI_{R_{\vec{v}}\T}\right]_{us} \right). \nonumber\\
\EEA

Let $\vec{x}=\vec{v}\amalg \vec{w}:=\vec{v}\cup\vec{w}\setminus\{i\in \vec{v}\cup\vec{w}\ |\ 
\exists j\in \vec{v}\cup\vec{w} \ |\ i\twoheadrightarrow j\}$. Then
one easily proves that
 $L_{\vec{w}}\circ L_{\vec{v}}(\T)=L_{\vec{x}}(\T)$,
$R_{\vec{v}}(\T)=R_{\vec{v}}\circ R_{\vec{x}}(\T)$  and $R_{\vec{w}}\circ L_{\vec{v}}(\T)=
L_{\vec{v}}\circ R_{\vec{x}}(\T)$. Hence
\BEA  [\del {\cal R} I_{\T}]_{t us}&=&
\sum_{\vec{x}\models V(\T)} [{\cal R} I_{L_{\vec{x}}\T}]_{t u}
\left( -\sum_{\vec{v}\models V(R_{\vec{x}}\T)}  [{\cal R}I_{L_{\vec{v}}(R_{\vec{x}}\T)}]_{us}
[{\cal R} \SkI_{R_{\vec{v}}(R_{\vec{x}}\T)}]_{s} + [\del {\cal R}\SkI_{R_{\vec{x}}\T}]_{us}\right) \nonumber\\
&=& \sum_{\vec{x}\models V(\T)} [{\cal R}I_{L_{\vec{x}}\T}]_{t u} [{\cal R} I_{R_{\vec{x}}\T}]_{us}.
\EEA

\hfill \eop


\subsection{Permutation graphs}


Consider now a permutation $\sigma\in\Sigma_n$. Applying Fubini's theorem yields
\BEA I_n^{\ell}(\Gamma)&=& \int_s^t d\Gamma_{x_1}(\ell(1))\int_s^{x_1} d\Gamma_{x_2}(\ell(2))
\ldots \int_s^{x_{n-1}} d\Gamma_{x_{n}}(\ell(n)) \nonumber\\
&=& \int_{s_1}^{t_1} d\Gamma_{x_{\sigma(1)}}(\ell(\sigma(1)))\int_{s_2}^{t_2} d\Gamma_{x_{\sigma(2)}}
(\ell(\sigma(2))) \ldots  \int_{s_{n}}^{t_{n}} d\Gamma_{x_{\sigma(n)}}(\ell(\sigma(n))), 
\nonumber\\ \label{eq:2:2.8}
\EEA
with $s_1=s$, $t_1=t$ and $s_j\in\{s\}\cup\{x_{\sigma(i)}, i<j\}$,
 $t_j\in\{t\}\cup\{x_{\sigma(i)},i<j\}$
$(j\ge 2)$. Now decompose $\int_{s_j}^{t_j} d\Gamma_{x_{\sigma(j)}}(\ell(\sigma(j)))$ into 
$\left( \int_s^{t_j}-\int_s^{s_j}\right) d\Gamma_{x_{\sigma(j)}}(\ell(\sigma(j)))$ if
$s_j\not=s,t_j\not=t$, and $\int_{s_j}^{t} d\Gamma_{x_{\sigma(j)}}(\ell(\sigma(j)))$ into 
$\left( \int_s^{t}-\int_s^{s_j}\right) d\Gamma_{x_{\sigma(j)}}(\ell(\sigma(j)))$ if
$s_j\not=s$. Then $I_n^{\ell}(\Gamma)$ has been rewritten as a sum of terms of the form
\BEQ\pm \int_s^{\tau_1} d\Gamma_{x_1}(\ell(\sigma(1)))\int_s^{\tau_2}d\Gamma_{x_2}(\ell(\sigma(2)))
\ldots\int_s^{\tau_{n}} d\Gamma_{x_{n}}(\ell(\sigma(n))), \label{eq:2:2.9} \EEQ
where $\tau_1=t$ and $\tau_j\in\{t\}\cup\{x_i, i<j\}$, $j=2,\ldots,n$. Note the renaming of 
variables and vertices from eq. (\ref{eq:2:2.8}) to eq. (\ref{eq:2:2.9}).  Encoding each of these
expressions by the forest  $\T$ with set of vertices  $V(\T)=\{1,\ldots,n\}$, label function
$\ell\circ\sigma$,  roots $\{j=1,\ldots,n\ |\
\tau_j=t\}$, and oriented edges $\{(j,j^-)\ |\ j=2,\ldots,n, \tau_j=x_{j^-}\}$, yields
\BEQ I_n^{\ell}(\Gamma)=I_{\T^{\sigma}}(\Gamma) \EEQ
for some $\T^{\sigma}\in{\cal T}$ called {\bf  permutation graph associated to $\sigma$}. 

Summarizing:

\begin{Lemma}[permutation graphs] \label{lem:2:sigma}

To every permutation $\sigma\in\Sigma_n$ is associated a permutation graph
\BEQ \T^{\sigma}=\sum_{j=1}^{J_{\sigma}} g(\sigma,j) \T_j^{\sigma}\in {\cal T},\EEQ
$g(\sigma,j)=\pm 1$, each forest $\T_j^{\sigma}$ being provided
by construction  with a total ordering compatible with its tree structure, image of the
ordering $\{v_1<\ldots<v_n\}$ of the trunk tree $\T^{\Id_n}$ by the permutation $\sigma$. The
label function of $\T^{\sigma}$ is $\ell\circ\sigma$, where $\ell$ is the original label
function of $\T^{\Id_n}$.

\end{Lemma}

\medskip

\begin{figure}[h]
  \centering
   \includegraphics[scale=0.35]{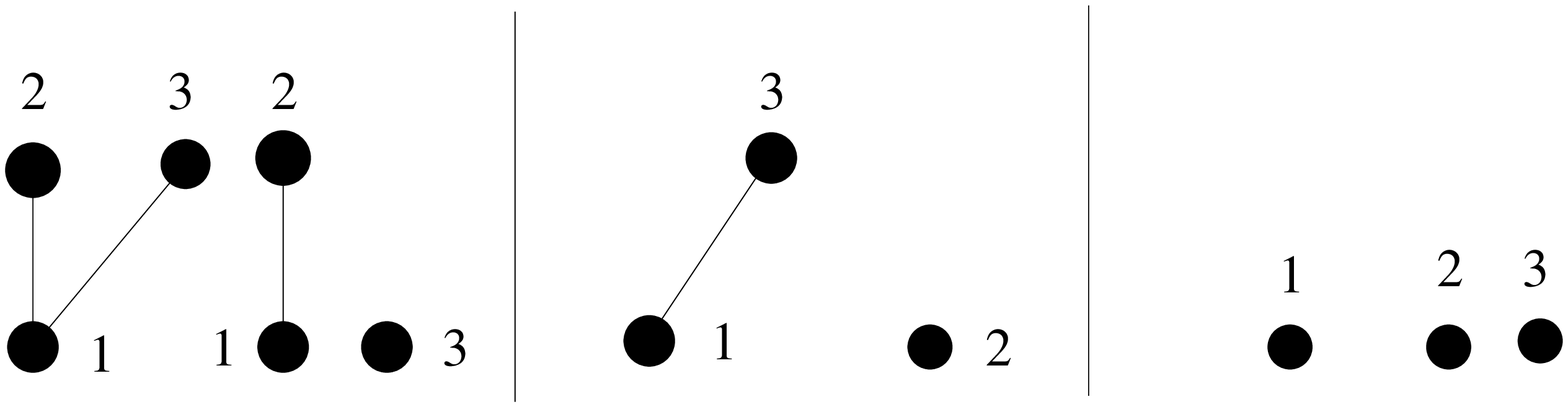}
   \caption{\small{Example \ref{ex:2:1}. From left to right: $\T^{\sigma}_1, \T^{\sigma}_2$;
 $L_{\{1\}}\T^{\sigma}_1 \otimes R_{\{1\}}\T^{\sigma}_1$; $L_{\{1,2\}}\T^{\sigma}_1
\otimes R_{\{1,2\}}\T^{\sigma}_1$ }}
 \label{Fig4bis}
\end{figure}

\begin{Example} \label{ex:2:1}

 Let $\sigma=\left(\begin{array}{ccc} 1 & 2 & 3\\ 2 & 3 & 1 \end{array}\right)$.
Then 
\BEA &&  \int_s^t d\Gamma_{x_1}(\ell(1)) \int_s^{t_2} d\Gamma_{x_2}(\ell(2)) \int_s^{t_3}
d\Gamma_{x_3}(\ell(3))= \nonumber\\
&& \qquad -\int_s^t d\Gamma_{x_2}(\ell(2)) \int_s^{x_2} d\Gamma_{x_3}(\ell(3))
\int_s^{x_2} d\Gamma_{x_1}(\ell(1)) \nonumber\\
&& \qquad \qquad  + \int_s^t d\Gamma_{x_2}(\ell(2)) \int_s^{x_2}
d\Gamma_{x_3}(\ell(3)) \ .\ \int_s^t d\Gamma_{x_1}(\ell(1)).\EEA
Hence $\T^{\sigma}=-\T_1^{\sigma}+\T_2^{\sigma}$ is the sum of a tree and of a forest with two
components (see Fig. \ref{Fig4bis}). 

\end{Example}

\subsection{Fourier normal ordering algorithm}


Let $\Gamma=(\Gamma(1),\ldots,\Gamma(d))$ be a compactly supported, 
smooth path, and ${\bf \Gamma}^n(i_1,\ldots,i_n)$
some iterated integral of $\Gamma$. To regularize ${\bf \Gamma}^n(i_1,\ldots,i_n)$, we shall
apply the following algorithm (a priori formal, since skeleton integrals may be infra-red divergent) :

\begin{enumerate}
\item (Fourier projections) Split the measure $\mu=d\Gamma(i_1)\otimes\ldots\otimes d\Gamma(i_n)$ into
$\sum_{\sigma\in\Sigma_n} {\cal F}^{-1} \left({\bf 1}_{D^{\sigma}}(\xi)\hat{\mu}(\xi) \right)$,
where $D^{\sigma}=\{(\xi_1,\ldots,\xi_n)\in\R^n \ |\ |\xi_{\sigma(1)}|
\le \ldots\le |\xi_{\sigma(n)}|\}$, and $\hat{\mu}$ is the Fourier transform of $\mu$. We shall
write 
\BEQ \mu^{\sigma}:={\cal F}^{-1} \left({\bf 1}_{D^{\sigma}}.\hat{\mu}\right)\circ\sigma=
{\cal F}^{-1} \left( {\bf 1}_{D^{\Id_n}} . (\hat{\mu}\circ\sigma)\right); \label{eq:2:mu-sigma} \EEQ

\item Rewrite $I_n^{\ell} \left({\cal F}^{-1} ({\bf 1}_{D^{\sigma}}.\hat{\mu}) \right)$, where
$\ell(j)=i_j$, as $I_{\T^{\sigma}}(\mu^{\sigma}):=\sum_{j=1}^{J_{\sigma}}
g(\sigma,j)I_{\T^{\sigma}_j}(\mu^{\sigma})$, where $\T^{\sigma}$ is the permutation
graph defined in subsection 2.4;

\item Replace $I_{\T^{\sigma}}(\mu^{\sigma})$ with some regularized integral as
in Definition \ref{def:2:reg} and  Lemma \ref{lem:2:reg},
\BEQ {\cal R}I_{\T^{\sigma}}(\mu^{\sigma}):=\sum_{j=1}^{J_{\sigma}} g(\sigma,j) 
{\cal R}I_{\T_j^{\sigma}}(\mu^{\sigma});\EEQ

\item Sum the terms corresponding to all possible permutations, yielding ultimately
\BEQ {\cal R}{\bf\Gamma}^n(i_1,\ldots,i_n)=\sum_{\sigma\in\Sigma_n} {\cal R}I_{\T^{\sigma}}(\mu^{\sigma}).\EEQ

\end{enumerate}

Explicit formulas for $\Gamma=B^{\eta}$ may be found in the following section.

\medskip

\begin{Theorem} \cite{Unt09bis}

${\cal R}{\bf\Gamma}$ satisfies the multiplicative (ii) and geometric (iii) properties defined
in the Introduction.

\end{Theorem}

The proof given in \cite{Unt09bis} shows actually that {\em any}
choice of linear  maps $[{\cal R}\Sk I_{\T}]_t:\mu\to [{\cal R}\Sk I_{\T}(\mu)]_t$
such that 

(i) $[{\cal R}\Sk I_{\T_1.\T_2}(\mu_1\otimes\mu_2)]_t=[{\cal R}\Sk I_{\T_1}(\mu_1)]_t [{\cal R}\Sk I_{\T_2}(\mu_2)]_t$
and

 (ii) $[{\cal R}\Sk I_{\T} (f)]_t=[\Sk I_{\T}(f)]_t=\int^t f(u)\ du$ if $\T$ is the trivial
tree with one vertex,

yields a regularized rough path over $\Gamma$ if $\Gamma$ is {\em smooth}. Hence our
 'cut' Fourier domain construction is arbitrary if convenient. As already said in the 
Introduction, it seems natural to look for some more restrictive rules for the regularization;
iterated renormalization schemes (such as BPHZ or dimensional regularization) are obvious
candidates (work in progress). The question is: is such or such regularization scheme better
in any sense ? Contrary to the case of quantum field theory where all renormalization
schemes may
be implemented by local counterterms, which amount to a change of the value of the
(finite number of) parameters
in the functional integral (which are experimentally measurable),
 and give ultimately after resumming the perturbation series
one and only one theory, we do not know of any {\em probabilistically motivated} reason
to choose a particular regularization scheme here.


\section{Rough path construction for fBm: case of distinct indices}


The strategy is now to choose an appropriate regularization procedure, so that regularized
skeleton
integrals of $B^{\eta}$ are finite and satisfy the  uniform H\"older and
convergence rate  estimates given in Theorem \ref{th:0}.


\subsection{Analytic approximation of fBm}


Recall $B$ may be defined via the harmonizable representation \cite{SamoTaq}
\BEQ B_t=c_{\alpha} \int_{\R}  |\xi|^{\half-\alpha}  \frac{e^{\II t\xi}-1}{\II \xi}
\ W(d\xi) \EEQ
where $(W_{\xi},\xi\ge 0)$ is a complex Brownian motion extendeded to $\R$
 by setting $W_{-\xi}=-\overline{W}_{\xi}$ $(\xi\ge 0)$, and 
$c_{\alpha}=\half \sqrt{-\frac{\alpha}{\cos\pi\alpha\Gamma(-2\alpha)}}$.

We shall use the following approximation of $B$ by a family of centered Gaussian
processes $(B^{\eta},\eta>0)$ living in the first chaos of $B$.

\begin{Definition}[approximation $B^{\eta}$]

Let, for $\eta>0$,

\BEQ B_t^{\eta}=c_{\alpha} \int_{\R} e^{-\eta|\xi|} |\xi|^{\half-\alpha} 
 \frac{e^{\II t\xi}-1}{\II \xi}\ W(d\xi). \EEQ

\end{Definition}

The process $B^{\eta}$ is easily seen to have a.s. smooth paths. The infinitesimal covariance
$\esper (B^{\eta})'_s (B^{\eta})'_t$ may be computed explicitly using the Fourier transform
\cite{Erd54}
\BEQ  {\cal F}K^{',-}_{\eta}(\xi)=\frac{1}{\sqrt{2\pi}} \int_{\R} K^{',-}_{\eta}(x) e^{-\II x\xi} dx
=-\frac{\pi\alpha}{2\cos\pi\alpha\Gamma(-2\alpha)} e^{-2\eta |\xi|} |\xi|^{1-2\alpha} 
{\bf 1}_{|\xi|>0},\EEQ 
where $K^{',-}_{\eta}(s-t):=\frac{\alpha(1-2\alpha)}{2\cos\pi\alpha} (-\II(s-t)+2\eta)^{2\alpha-2}$.
By taking the real part of these expressions, one finds that $B^{\eta}$ has the same law as
the analytic approximation of $B$ defined in \cite{Unt08}, namely, $B^{\eta}=
\Gamma_{t+\II\eta}+\Gamma_{t-\II\eta}=2\Re\Gamma_{t+\II\eta}$, where $\Gamma$ is the
analytic fractional Brownian motion (see also \cite{TinUnt08}).

\medskip


\subsection{Choice of the regularization procedure}


 Let $\sigma\in\Sigma_n$ be a permutation. Recall (see Lemma \ref{lem:2:sigma}) that
the permutation graph  $\T^{\sigma}$
 may be written as  a finite sum
$\sum_{j=1}^{J_{\sigma}} g(\sigma,j) \T^{\sigma}_j$, where each $\T^{\sigma}_j$ is a forest which
is automatically provided with a total ordering. In the two following  subsections, 
we shall consider regularized tree or skeleton integrals, ${\cal R}I_{\T}$ or 
${\cal R}\Sk I_{\T}$, for a forest $\T$ which is one of the $\T^{\sigma}_j$.

\begin{Definition}  \label{def:6:RTreg}
{\it Fix $C_{reg}\in(0,1).$ Let, for  $\T$ with set of vertices $V(\T)=\{v_1<\ldots<v_j\}$,
\BEQ \R_{+}^{\T}:=\big\{(\xi_{v_1},\ldots,\xi_{v_j})\in\R^{\T}\ |\ |\xi_{v_1}|\le\ldots\le|\xi_{v_j}|\},\EEQ
\BEQ \R_{reg}^{\T}:=\big\{(\xi_{v_1},\ldots,\xi_{v_j})\in\R^{\T}_+\ | \  \forall v\in V(\T),
|\xi_v+\sum_{w\twoheadrightarrow v} \xi_w|>C_{reg}\max\{|\xi_w|;\ w\twoheadrightarrow v\}
 \ \big\}, \label{eq:3:RTreg} \EEQ
and ${\cal R}I_{\T}$,  resp.  ${\cal R}\SkI_{\T}$ be the corresponding 
$\R^{\T}_{reg}$-regularized iterated, resp.
skeleton integrals as in subsection 2.3.
}
\end{Definition}

Condition (\ref{eq:3:RTreg}) ensures that the denominators in the skeleton integrals are
not too small (see Lemma \ref{lem:2:SkI}).

\medskip

The following Lemma (close to arguments used in the study of random Fourier series \cite{Kah})
is fundamental for the estimates of the following subsections.

\begin{Lemma} \label{lem:3:Kah}

\begin{itemize}

\item[(i)] Let $F(u)=\int_{\R} dW_{\xi} a(\xi)e^{\II u\xi}$, where $|a(\xi)|^2\le C|\xi|^{-1-2\beta}$
for some $0<\beta<1$: then, for every $u_1,u_2\in\R$, 
\BEQ \esper |F(u_1)-F(u_2)|^2 \le C' |u_1-u_2|^{2\beta}.\EEQ

\item[(ii)] Let $\tilde{F}(\eta)=\int_{\R} dW_{\xi} a(\xi)e^{-\eta|\xi|}$ $(\eta>0)$,
 where $|a(\xi)|^2\le C|\xi|^{-1-2\beta}$
for some $0<\beta<1$: then, for every $\eta_1,\eta_2\in\R_+$, 
\BEQ \esper |\tilde{F}(\eta_1)-\tilde{F}(\eta_2)|^2 \le C' |\eta_1-\eta_2|^{2\beta}.\EEQ

\end{itemize}

\end{Lemma}

{\bf Proof.} Bound $|e^{\II u_1\xi}-e^{\II u_2\xi}|$ by $|u_1-u_2| |\xi|$ for $|\xi|\le\frac{1}{|u_1-u_2|}$ and by $2$ otherwise, and similarly for $|e^{-\eta_1|\xi|}-e^{-\eta_2|\xi|}|$.
 Note the variance integral is infra-red convergent near $\xi=0$. 
\hfill \eop

{\bf Remark:} Unless $|a(\xi)|^2$ is $L^1_{loc}$ near $\xi=0$, only the {\em increments} $F(u_1)-F(u_2)$,
$\tilde{F}(\eta_1)-\tilde{F}(\eta_2)$ are well-defined.


\subsection{Estimates for the increment term}


In this paragraph, as in the next one, we consider regularized tree integrals associated
to ${\cal R}{\bf B}^{n,\eta}(i_1,\ldots,i_n)$ where $i_1\not=\ldots\not=i_n$ are distinct indices,
so that $B(i_1),\ldots,B(i_n)$ are {\it independent}.

\begin{Lemma}[H\"older estimate and rate of convergence] \label{lem:3:increment-Holder-rate}
{\it Let $\T=\T^{\sigma}_j$ for some $j$,  and $\alpha<1/|V(\T)|$.
\begin{enumerate}
\item
 The skeleton term  
\BEQ [G^{\eta,\sigma}_{\T}(i_1,\ldots,i_n)]_u=
\left[{\cal R}\SkI_{\T}\left( 
\left( dB^{\eta}(i_1)\otimes\ldots\otimes dB^{\eta}(i_n) \right)^{\sigma}
  \right) \right]_u  \EEQ  (see eq. (\ref{eq:2:mu-sigma}))   writes
\BEA && [G^{\eta,\sigma}_{\T}(i_1,\ldots,i_n)]_u=(-\II c_{\alpha})^{|V(\T)|} 
\int\ldots\int_{(\xi_v)_{v\in V(\T)}\in \R^{\T}_{reg}} \prod_{v\in V(\T)} dW_{\xi_v}(i_{\sigma(v)}) \nonumber\\
&& \qquad \qquad \qquad  e^{\II u \sum_{v\in V(\T)} \xi_v} e^{-\eta\sum_{v\in V(\T)} |\xi_v|}
\frac{\prod_{v\in V(\T)} |\xi_v|^{\half-\alpha}}
{\prod_{v\in V(\T)} \left[ \xi_v+\sum_{w\twoheadrightarrow v} \xi_w\right]}. \label{eq:6:40}\EEA
\item
It satisfies the uniform H\"older estimate:
\BEQ \esper \left| [\del G^{\eta,\sigma}_{\T}(i_1,\ldots,i_n)]_{ts} \right|^2
    \le C|t-s|^{2\alpha|V(\T)|}.\EEQ
\item (rate of convergence) :  there exists a constant $C>0$ such that, for every $\eta_1,
\eta_2>0$ and $s,t\in\R$,
\BEQ \esper \left|  [\del G^{\eta_1,\sigma}_{\T}(i_1,\ldots,i_n)]_{ts}-
 [\del G^{\eta_2,\sigma}_{\T}(i_1,\ldots,i_n)]_{ts}
 \right|^2 \le C|\eta_1-\eta_2|^{2\alpha}.\EEQ
\end{enumerate}
}
\end{Lemma}

{\bf Proof.}  

\begin{enumerate}

\item Follows from  Lemma \ref{lem:2:SkI} and the definitions of $B^{\eta}$ and of regularized
integrals in the previous subsections  2.3 and 3.1.

\item (H\"older estimate) 

One may just as well (by multiplying the integral estimates on each tree component) assume $\T$
is a tree, i.e. $\T$ is connected.

Let $V(\T)=\{v_1<\ldots<v_{|V(\T)|}\}$, so that $|\xi_{v_1}|\le \ldots\le |\xi_{v_{|V(\T)|}}|$.
  Since every vertex $v\in V(\T)\setminus\{v_1\}$ connects
to the root $v_1$, one has 
\BEQ |V(\T)|\ .\ |\xi_{v_{|V(\T)|}}|\ge |\xi_{v_1}+\ldots+\xi_{v_{|V(\T)|}}|>C_{reg}|\xi_{v_{|V(\T)|}}|,
\EEQ so
that $\xi:=\sum_{v\in V(\T)} \xi_v$ 
 is comparable to $\xi_{v_{|V(\T)|}}$, i.e. belongs to $[C^{-1}\xi_{v_{|V(\T)|}},C\xi_{v_{|V(\T)|}}]$ if
$C$ is some large enough positive constant. Write $[G^{\eta,\sigma}_{\T}(i_1,\ldots,i_n)]_u
=\int_{\R} e^{\II u\xi}
a(\xi) d\xi$.

Vertices at which 2 or more branches join are called {\it nodes}, and vertices to which no vertex
is connected are called {\it leaves} (see Fig. \ref{Fig5}).

\begin{figure}[h]
  \centering
   \includegraphics[scale=0.35]{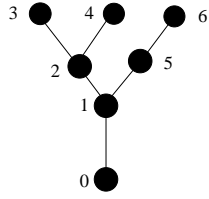}
   \caption{\small{3,4,6 are leaves; 1, 2 and 5 are nodes, 2 and 5 are uppermost; branches are e.g. 
$Br(2\twoheadrightarrow 1)$ or $Br(6\twoheadrightarrow 1)$.}}
  \label{Fig5}
\end{figure}

 The set $Br(v_1\twoheadrightarrow v_2)$ of vertices from a leaf
or a node $v_1$ to a node $v_2$ (or to the root) is called a {\it branch}
 if it does not contain any other node. By
convention, $Br(v_1\twoheadrightarrow v_2)$ includes $v_1$ and excludes $v_2$.

Consider an uppermost node $n$, i.e. a node to which no other node is connected, together with the set
of leaves $\{w_1<\ldots<w_J\}$ above $n$. Let $p_j=|V(Br(w_j\twoheadrightarrow n))|$. Note
that $\left( \frac{|\xi_{n}|^{\half-\alpha}}{\xi_n+\sum_{w\twoheadrightarrow n} \xi_w} \right)^2\lesssim 
|\xi_{w_J}|^{-1-2\alpha}$. Now we proceed to estimate $\Var\  a(\xi)$.
On the branch number $j$ from $w_j$ to $n$,
\BEA &&  \int\ldots\int_{|\xi_v|\le |\xi_{w_j}|, v\in Br(w_j\twoheadrightarrow n)\setminus\{w_j\} }
\left[ \prod_{v\in Br(w_j\twoheadrightarrow n)} \frac{e^{-\eta |\xi_v|} |\xi_v|^{\half-\alpha}}{\xi_v+\sum_{w\twoheadrightarrow v}
\xi_w} \right]^2  \nonumber\\
&& \qquad \lesssim |\xi_{w_j}|^{-1-2\alpha p_j} \EEA
and (summing over  $\xi_{w_1},\ldots,\xi_{w_{J-1}}$  and  over $\xi_n$)
\BEA && |\xi_{w_J}|^{-1-2\alpha p_J} \int_{|\xi_{w_{J-1}}| \le |\xi_{w_J}|} d\xi_{w_{J-1}}
|\xi_{w_{J-1}}|^{-1-2\alpha p_{J-1}} \ \nonumber\\
&& \qquad \left( \ldots
\left( \int_{|\xi_{w_1}|\le|\xi_{w_2}|} d\xi_{w_1} |\xi_{w_1}|^{-1-2\alpha p_1} \left(
 \int_{|\xi_n| \le |\xi_{w_1}|}
d\xi_n \frac{|\xi_n|^{1-2\alpha}}{\xi_{w_J}^2} \right)\right) \ldots \right) \nonumber\\
&& \lesssim |\xi_{w_J}|^{-(1+2\alpha p_j)+[2-2\alpha(1+p_1+\ldots+p_{J-1})]-2}
=|\xi_{w_J}|^{-1-2\alpha W(n)},\EEA
where $W(n)=p_1+\ldots+p_J+1=|\{v: v\twoheadrightarrow n\}|+1$ is the {\it weight} of $n$.

One may then consider the reduced tree $\T_n$ obtained by shrinking all vertices above $n$ (including
$n$) to {\it one} vertex  with weight $W(n)$ and perform the same operations on $\T_n$. Repeat this
inductively until $\T$ is shrunk to one point. In the end, one gets
$\Var \ a(\xi) \lesssim |\xi_{v_{|V(\T)|}}|^{-1-2\alpha|V(\T)|} \lesssim |\xi|^{-1-2\alpha|V(\T)|}$. 
Now apply
Lemma \ref{lem:3:Kah} (i).

\bigskip

\item (rate of convergence)

Let $X_u^{\eta_1,\eta_2}:=[G^{\eta_1,\sigma}_{\T}(i_1,\ldots,i_n)]_{u}-
 [ G^{\eta_2,\sigma}_{\T}(i_1,\ldots,i_n)]_{u}$. Expanding $\prod_{j=1}^{|V(\T)|}
e^{-\eta_1 |\xi_{j}|} -\prod_{j=1}^{|V(\T)|} e^{-\eta_2 |\xi_j|}$ as
\BEQ \sum_{j=1}^{|V(\T)|} e^{-\eta_2(|\xi_{v_1}|+\ldots+|\xi_{v_{j-1}}|)}
(e^{-\eta_1|\xi_{v_j}|}-e^{-\eta_2|\xi_{v_j}|}) e^{-\eta_1(|\xi_{v_{j+1}}|+\ldots
+|\xi_{v_{V(\T)}}|)} \nonumber\\ \EEQ
 gives $X_u^{\eta_1,\eta_2}$ as a sum, $X_u^{\eta_1,\eta_2}=\sum_{v\in V(\T)} 
X_u^{\eta_1,\eta_2}(v)$, where $X_u^{\eta_1,\eta_2}(v)=\int d\xi_v b_u(\xi_v)
(e^{-\eta_1 |\xi_v|}-e^{-\eta_2|\xi_v|}) $ 
is obtained from $[G^{\eta,\sigma}_{\T}(i_1,\ldots,i_n)]_{u}$
by replacing $e^{-\eta|\xi_v|}$ with $e^{-\eta_1|\xi_v|}-e^{-\eta_2|\xi_v|}$,
 and $e^{-\eta|\xi_w|}$, $w\not=v$ either by
$e^{-\eta_1|\xi_w|}$ or by $e^{-\eta_2|\xi_w|}$. We want to estimate $\Var\  b_u(\xi_v)$
 uniformly in $u$.

\medskip 

Fix the value of $\xi_v$ in the  computations in the above  proof for the H\"older estimate.
Let $w_J$ be the maximal leaf above $v$, and $n\twoheadrightarrow v$ be the node just above $v$ if $v$ is not a node,
$n=v$ otherwise.
 Summing over all
nodes above $n$ and taking the variance leads to an expression bounded by $|\xi_{w_J}|^{-1-2\alpha W(n)}$,
where  $W(n)=|\{w\ :\ w\twoheadrightarrow n\}|+1$ is as before
  the weight of $n$. Consider now the corresponding shrunk
tree $\T_n$.   Let $\T_n(v)$  be the trunk tree defined by
$\T_n(v)=\{w\in\T_n: w\twoheadrightarrow v \ {\mathrm{or}}\ v\twoheadrightarrow w\}\cup\{v\}$; similarly,
let $\T(v)$ be the tree defined by $\T(v)=\{w\in\T: w\twoheadrightarrow v \ {\mathrm{or}}\ v\twoheadrightarrow w\}\cup\{v\}$, so that $\T_n(v)$ is the corresponding shrunk tree.
Sum over all vertices $w\in \T_n(v)\setminus\{v\}$. The variance of the coefficient of $e^{-\eta_1 |\xi_v|}$ is
\BEA S(\xi_v) &\lesssim& \int_{|\xi_n|\ge |\xi_{v}|} d\xi_n |\xi_n|^{-1-2\alpha W(n)}|\xi_n|^{-1-2\alpha}
 \int_{|\xi_w|\le |\xi_n|,w\in \T_n(v)\setminus\{n,v\}}  \nonumber\\
&& \qquad \left[
\prod_{w\in\T_n(v)\setminus\{n,v\}} d\xi_w \ .\  |\xi_n|^{-(1+2\alpha)} \right]  \nonumber\\
&\lesssim &  \int_{|\xi_n|\ge |\xi_{v}|} d\xi_n |\xi_n|^{-2-2\alpha|\T(v)|} \lesssim |\xi_v|^{-1-2\alpha|\T(v)|}  \EEA
if $v\not =n$, and
\BEQ S(\xi_v)\lesssim |\xi_n|^{-1-2\alpha W(n)} 
\int_{|\xi_w|\le |\xi_n|,w\in\T_n(v)\setminus\{n\}} \prod_{w\in \T_n(v)\setminus\{n\}}
|\xi_n|^{-(1+2\alpha)} \lesssim |\xi_v|^{-1-2\alpha|\T(v)|} \EEQ
if $v=n$.

 Removing
the vertices belonging to $\T(v)$ from $\T$ leads to a forest which gives a finite contribution
to the variance. Hence (by Lemma \ref{lem:3:Kah} (ii)) $\esper|X_u^{\eta_1,\eta_2}(v)|^2\lesssim 
|\eta_1-\eta_2|^{2\alpha|\T(v)|}.$

\end{enumerate}

 \hfill \eop

The notion of {\em weight} $W(v)$ {\em of a vertex} $v$ introduced in this proof will be used again in
subsections 3.4 and  4.1.


\subsection{Estimates for boundary terms}


Let $\T=\T^{\sigma}_j$ for some $\sigma\in \Sigma_n$, and $i_1\not=\ldots\not=i_n$
as in the previous subsection.  By multiplying the estimates on each tree component,
 one may just as well assume $\T$ is a tree, i.e. is connected.

We shall now prove estimates for the boundary term ${\cal R}I_{\T} \left( \left( dB^{\eta}(i_1)\otimes \ldots\otimes dB^{\eta}(i_n) \right)^{\sigma} 
\right) (\partial)$
associated to $\T$ (see Lemma \ref{lem:2:reg}).

\begin{Lemma}

{\it Let $\T=\T^{\sigma}_j$ for some $j$ (so that $n=|V(\T)|$).
\begin{enumerate}

\item (H\"older estimate)
 The regularized
 boundary term ${\cal R}I_{\T} \left( \left( dB^{\eta}(i_1)\otimes \ldots\otimes dB^{\eta}(i_n) \right)^{\sigma} 
\right)(\partial)$ satisfies:
 \BEQ \esper \left|[ {\cal R} I_{\T}\left( \left( dB^{\eta}(i_1)\otimes \ldots\otimes dB^{\eta}(i_n) \right)^{\sigma} 
\right)(\partial) \right]_{ts}|^2
\le C |t-s|^{2\alpha|V(\T)|} \label{eq:6:boundary-Holder} \EEQ
for a certain constant $C$.

\item (rate of convergence)
 There exists a positive constant $C$ such that, for every $\eta_1,\eta_2>0$,
\BEA &&  \esper \left|  [{\cal R}I_{\T}\left( \left( dB^{\eta_1}(i_1)\otimes \ldots\otimes dB^{\eta_1}(i_n) \right)^{\sigma} 
\right)(\partial)]_{ts} - \right.\nonumber\\  && \left. \qquad \qquad 
-[{\cal R}I_{\T}\left( \left( dB^{\eta_2}(i_1)\otimes \ldots\otimes dB^{\eta_2}(i_n)
 \right)^{\sigma} 
\right)(\partial)]_{ts} \right|^2 
 \le C|\eta_1-\eta_2|^{2\alpha}. \nonumber\\ \EEA
\end{enumerate}

\label{lem:3:boundary-Holder-rate}
}

\end{Lemma}

{\bf Proof.}

\begin{enumerate}

\item

 Apply repeatedly Lemma \ref{lem:2:reg} to $\T$: in the end,
 $[{\cal R}I_{\T}\left( \left( dB^{\eta}(i_1)\otimes \ldots\otimes dB^{\eta}(i_n) \right)^{\sigma} 
\right)(\partial)]_{ts}$ appears as a sum of 'skeleton-type' terms of the form
(see Figure \ref{Fig6})
\BEA &&   A_{ts}:=[\del{\cal R}\SkI_{ L\T}]_{ts}\ . \
 [{\cal R}\SkI_{R_{\vec{v}_l}\circ L_{\vec{v}_{l-1}}\circ\ldots\circ L_{\vec{v}_1}(\T)}]_{s}
\ldots
 [{\cal R}\SkI_{R_{\vec{v}_2}\circ L_{\vec{v}_1}(\T)}]_{s}
[{\cal R}\SkI_{ R_{\vec{v}_1\T}}]_{s} \nonumber\\
&& \qquad \qquad \qquad \qquad
\left( \left( dB^{\eta}(i_1)\otimes \ldots\otimes dB^{\eta}(i_n) \right)^{\sigma} 
\right), \nonumber\\ \label{eq:6:skeleton-type}
\EEA
where $\vec{v}_1=(v_{1,1}<\ldots<v_{1,J_1})\models \T$, $\vec{v}_2\models L_{\vec{v}_1}\T$, $\ldots$,
$\vec{v}_l=(v_{l,1}<\ldots<v_{l,J_l})\models L_{\vec{v}_{l-1}}\circ \ldots \circ L_{\vec{v}_1}(\T)$ and
$L\T:=L_{\vec{v}_l}\circ\ldots\circ L_{\vec{v}_1}(\T)$. In eq. (\ref{eq:6:skeleton-type}) the forest $\T$
has been split into a number of sub-forests, $L\T\cup \left( \cup_{j=1}^J \T_j \right)$;
 we call this splitting {\em the splitting associated to} $A_{ts}$ for further reference.

\medskip

\underline{First step.}

 Let $ [B^{\vec{v}_1}_s [\vec{\xi}]]_u \prod_{j=1}^{J_1} dW_{\xi_{v_1,j}}
(i_{\sigma(v_{1,j})})$ 
be the contribution to    ${\cal R}\SkI_{R_{\vec{v}_1}\T}$ of all Fourier components such that
 $\vec{\xi}=(\xi_{v_{1,1}},\ldots,\xi_{v_{1,J_1}})$,  $|\xi_{v_{1,1}}|\le \ldots \le|\xi_{v_{1,J_1}}|$ is fixed. For definiteness (see Definition \ref{def:6:RTreg}),
\BEA && \left[{\cal R}\SkI_{R_{\vec{v}_1\T}} \left( \left( dB^{\eta}(i_1)\otimes \ldots\otimes dB^{\eta}(i_n) \right)^{\sigma} 
\right) \right]_u ((x_v)_{v\in L_{\vec{v}_1}\T}) \nonumber\\
 && = \int\ldots\int  {\bf 1}_{\R^{V(L_{\vec{v}_1}\T)\cup\vec{v}_1}_+}
\left((\xi_v)_{v\in V(L_{\vec{v}_1}\T)\cup\vec{v}_1}\right)
\left[ B_s^{\vec{v}_1}[\vec{\xi}] \prod_{j=1}^{J_1} dW_{\xi_{v_1,j}}
(i_{\sigma(v_{1,j})}) \right]  \ . \nonumber\\
&& \qquad .\  \left[ \prod_{v\in L_{\vec{v}_1}\T} c_{\alpha}
e^{-\eta |\xi_v|} e^{\II x_v\xi_v} |\xi_v|^{\half-\alpha} dW_{\xi_v}(i_{\sigma(v)}) \right].\EEA
 Then
\BEQ  \Var [B_s^{\vec{v}_1}[\vec{\xi}]]_{s} \lesssim 
\int\ldots\int \prod_{v\in\vec{v}_1} d\xi_v \left[ |\xi_v|^{-1-2\alpha}  \int\ldots\int_{|\xi_w|\ge|\xi_v|,w\in R_v\T\setminus\{v\}} \prod_{w\in R_v\T\setminus\{v\}} |\xi_w|^{-1-2\alpha} \right],
 \EEQ hence
\BEQ  \Var [B_s^{\vec{v}_1}[\vec{\xi}]]_{s} \lesssim 
\prod_{v\in\vec{v}_1}
|\xi_v|^{-2|V(R_v\T)|\alpha-1}.\EEQ

\medskip

\underline{Second step.}

More generally, let $ B_{s}^{\vec{v}_1,\ldots,\vec{v}_l}[\vec{\xi}]\prod_{j=1}^{J_l} dW_{\xi_{v_{l,j}}}(i_{\sigma(v_{l,j})})$ be the contribution to
\BEQ  [{\cal R}\SkI_{R_{\vec{v}_l}\circ L_{\vec{v}_{l-1}}\circ\ldots\circ L_{\vec{v}_1}(\T)}]_{s}
\ldots
 [{\cal R}\SkI_{R_{\vec{v}_2}\circ L_{\vec{v}_1}(\T)}]_{s}
[{\cal R}\SkI_{ R_{\vec{v}_1\T}}]_{s}
\left( \left( dB^{\eta}(i_1)\otimes \ldots\otimes dB^{\eta}(i_n) \right)^{\sigma} 
\right) \EEQ
  of all
Fourier components such that  $\vec{\xi}=(\xi_{v_{l,1}},\ldots,\xi_{v_{l,J_l}})$ is fixed.  Then
\BEQ \Var ( B^{\vec{v}_1,\ldots,\vec{v}_l}_{s}[\vec{\xi}]) \lesssim \prod_{v\in\vec{v}_l} |\xi_v|^{-2|V(R_{v}\T)|\alpha-1}\EEQ
(proof by induction on $l$).

\medskip

\underline{Third step.}

Let $V(L\T)=\{w_1<\ldots<w_{max}\}$. 
By definition, $A_{ts}=\int_{\R}  a_s(\Xi)(e^{\II \Xi t}-e^{\II \Xi s})d\Xi$,
with
\BEA &&  a_s(\Xi)=\int d\vec{\xi} \int\ldots\int_{( (\xi_w)_{w\in V(L\T)})\in D_{\vec{\xi}}} \prod_{w\in V(L\T)} dW_{\xi_w}(i_{\sigma(w)}) \nonumber\\
&& \qquad \qquad 
\frac{\prod_{w\in V(L\T)} (-\II c_{\alpha})e^{-\eta|\xi_w|} |\xi_w|^{\half-\alpha}}{\prod_{w\in V(L\T)}
(\xi_w+\sum_{w'\twoheadrightarrow w,w'\in V(L\T)} \xi_{w'})}  B_s^{\vec{v}_1,\ldots,
\vec{v}_l}[\vec{\xi}] \nonumber\\
\EEA
where Fourier components  in $D_{\vec{\xi}}$ satisfy in particular the following conditions:

\begin{itemize}
\item $|\xi_w+\sum_{w'\twoheadrightarrow w,w'\in V(L\T)} \xi_{w'}|>C_{reg} \max\{ |\xi_{w'}|: w'\twoheadrightarrow w, w'\in V(L\T)\}$;
in particular, $\left(\frac{|\xi_w|^{\half-\alpha}}{\xi_w+\sum_{w'\twoheadrightarrow w,w'\in
V(L\T)} \xi_{w'}}\right)^2 \lesssim |\xi_w|^{-1-2\alpha}$;
\item $\sum_{w\in V(L\T)} \xi_w=\Xi$;
\item for every $w\in V(L\T)$, $|\xi_w|\le |\xi_{w_{max}}|$ and $|\xi_w| \le |\xi_v|$ for every $v\in R(w):=\{
v=v_{l,1},\ldots,v_{l,J_l}\ |\ v\to w\}$ (note that $R(w)$ may be  empty). See Fig. \ref{Fig6}.
\end{itemize}

\begin{figure}[h]
  \centering
   \includegraphics[scale=0.35]{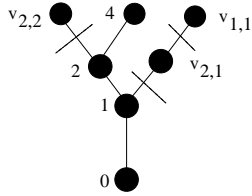}
   \caption{\small{Here $V(L\T)=\{0,1,2,4\}$, $R(0)=R(4)=\emptyset$, $R(1)=\{v_{2,1}\}$,
$R(2)=\{v_{2,2}\}$.}}
  \label{Fig6}
\end{figure}

Note that $|\Xi|\lesssim |\xi_{w_{max}}| \lesssim |\Xi|$ since every vertex in $V(L\T)$ connects to the root
(see first lines of the proof of Lemma \ref{lem:3:increment-Holder-rate} (2)).

If $w\in L\T$, split $R(w)$ into $R(w)_{>}\cup R(w)_{<}$, where $R(w)_{\gtrless}
:=\{v\in R(w)\ |\ v\gtrless w_{max}\}$. Summing over indices corresponding to vertices in $R\T_{>}:=
\{v=v_{l,1},\ldots,v_{l,J_l}\ |\ v>w_{max}\}=\cup_{w\in L\T}
R(w)_{>}$, one gets (see again proof of Lemma \ref{lem:3:increment-Holder-rate} (2))
\BEQ \prod_{v\in R\T_>} \int_{|\xi_v|\ge |\Xi|} d\xi_v |\xi_v|^{-2|V(R_v\T)|\alpha-1} \lesssim |\Xi|^{-2\alpha\sum_{v\in R\T_>} |V(R_v\T)|}.\EEQ

Let $w\in L\T\setminus\{w_{max}\}$ such that $R(w)_<\not=\emptyset$  (note that $R(w_{max})_<=\emptyset$).
 Let 
 $R(w)_< = \{v_{i_1}<\ldots<v_{i_j}\}$. Then (integrating over $(\xi_v), v\in R(w)_<$)
\BEA  && |\xi_w|^{-1-2\alpha} \int_{|\xi_{v_{i_1}}|\ge |\xi_w|} d\xi_{v_{i_1}} \int_{|\xi_{v_{i_2}}|\ge|\xi_{v_{i_1}}|} d\xi_{v_{i_2}}
 \ldots \int_{|\xi_{v_{i_j}}|\ge|\xi_{v_{i_{j-1}}}|} d\xi_{v_{i_j}} \nonumber\\
&& 
|\xi_{v_{i_1}}|^{-2|V(R_{v_{i_1}}\T)|\alpha-1} \ldots |\xi_{v_{i_j}}|^{-2|V(R_{v_{i_j}}\T)|\alpha-1} 
 \lesssim |\xi_w|^{-1-2\alpha(1+ \sum_{v\in R(w)_<} |V(R_v\T)|)}. \nonumber\\ \EEA

In other words, each vertex $w\in L\T$ 'behaves' as if it had a weight $1+\sum_{v\in R(w)_<} |V(R_v\T)|.$
Hence (by the same method as in the proof of Lemma \ref{lem:3:increment-Holder-rate} (2))
 $\Var(a_s(\xi))\lesssim |\Xi|^{-1-2\alpha(|V(L\T)|+\sum_{v\in R\T_<} |V(R_v\T)|)} \ .\ |\Xi|^{-2\alpha\sum_{v\in R\T_>} |V(R_v\T)|}
= |\Xi|^{-1-2\alpha|V(\T)|}.$ Now apply Lemma \ref{lem:3:Kah} (i).

\item  Similar to the 
proof of Lemma \ref{lem:3:increment-Holder-rate} (3). Details are left to the reader.
\end{enumerate}

\hfill \eop


\section{End of proof and final remarks}



\subsection{Estimates: case of coinciding indices}


Our previous estimates for $\esper|{\cal R}{\bf B}^{n,\eta}_{ts}(i_1,\ldots,i_{n})|^2$ 
(H\"older estimate) and \\
$\esper|{\cal R}{\bf B}^{n,\eta_1}_{ts}(i_1,\ldots,i_{n})-{\cal R}{\bf B}^{n,\eta_2}_{ts}(i_1,\ldots,i_{n})|^2$ (rate of convergence) with $i_1\not=\ldots\not=i_n$ rest on the independence of
the Brownian motions $W(i_1),\ldots,W(i_n)$.
 We claim that the same estimates also hold true for $\esper
|{\cal R}{\bf B}^{n,\eta}(i_1,\ldots,i_n)|^2$ and $\esper|{\cal R}{\bf B}^{n,\eta_1}_{ts}(i_1,\ldots,i_n)-{\cal R}{\bf B}^{n,\eta_2}_{ts}(i_1,\ldots,i_n)|^2$ if some of the indices
$(i_1,\ldots,i_n)$ coincide, with the {\it same definition} of the regularization procedure ${\cal R}$. The
key Lemma for the proof is

\begin{Lemma}[Wick's lemma](see \cite{LeB}, \S 5.1.2 and 9.3.4)

{\it Let $(X_1,\ldots,X_{n})$ be a centered Gaussian vector. Denote by $X_{i_1}\diamond\ldots\diamond X_{i_k}$
$(1\le i_1,\ldots,i_k\le n)$ or $:X_{i_1}\ldots X_{i_k}:$
 the Wick product of $X_{i_1},\ldots,X_{i_k}$ (also called: {\it normal ordering} of the product
$X_{i_1}\ldots X_{i_k}$), i.e. the projection
of the product $X_{i_1}\ldots X_{i_k}$ onto the $k$-th chaos of the Gaussian space generated by
$X_1,\ldots,X_{n}$. Then:

\begin{enumerate}

\item

\BEA
&& X_1\ldots X_{n} = X_{1}\diamond\ldots\diamond X_{n} 
  +  \sum_{(i_1,i_2)} \esper[X_{i_1}X_{i_2}]
X_{1}\diamond \ldots \diamond \check{X}_{i_1}\diamond \ldots\diamond \check{X}_{i_2}\diamond\ldots\diamond X_n \nonumber\\
&& \qquad  +  \ldots+ \sum_{(i_1,i_2),\ldots,(i_{2k+1},i_{2k+2})}
 \esper[X_{i_1}X_{i_2}]\ldots \esper[X_{i_{2k+1}}X_{i_{2k+2}}]  \nonumber\\
&& \qquad 
X_1\diamond\ldots\diamond \check{X}_{i_1}\diamond\ldots\diamond \check{X}_{i_2}
\diamond\ldots\diamond \check{X}_{i_{2k+1}} \diamond\ldots\diamond \check{X}_{i_{2k+2}}
\diamond\ldots\diamond X_n \nonumber\\
&& +\ldots, \EEA
where the sum ranges over all partial  pairings of indices $(i_1,i_2),\ldots,(i_{2k+1},i_{2k+2})$
$(1\le k\le \lfloor \frac{n}{2}\rfloor -1)$.

\item
For every set of indices $i_1,\ldots,i_{j},i'_1,\ldots,i'_{j}$,
\BEQ \esper\left[ (X_{i_1}\diamond\ldots\diamond X_{i_{j}})(X_{i'}\diamond\ldots\diamond X_{i'_{j}}) \right]
=\sum_{\sigma\in \Sigma_j} \prod_{m=1}^{j} \esper[ X_{i_m} X_{i'_{\sigma(m)}}]. \EEQ
\end{enumerate}
}
\label{lem:7:Wick}
\end{Lemma}

In our case (considering ${\cal R} {\bf B}^{n,\eta}_{ts}(i_1,\ldots,i_n)$) we get a decomposition of the product $dW_{\xi_1}(i_1)\ldots dW_{\xi_{n}}(i_n)$
into $dW_{\xi_1}(i_1)\diamond\ldots\diamond dW_{\xi_n}(i_n)$,
 plus the sum over all possible non-trivial
pair contractions, schematically  $\langle W'_{\xi_j}(i_j)  W'_{\xi_{j'}}(i_{j'})\rangle=\del_0(\xi_j+\xi_{j'})
\del_{i_j,i_{j'}}.$

Consider first the normal ordering of ${\cal R}{\bf B}^{n,\eta}_{ts}(i_1,\ldots,i_n)$. As in
the proof of Lemma 5.10 in \cite{TinUnt08}, let $\Sigma_{\vec{i}}$ be the
 'index-fixing' subgroup of 
$\Sigma_n$ such that : $\sigma'\in\Sigma_{\vec{i}}\Longleftrightarrow \forall j=1,\ldots,n,
\ i_{\sigma'(j)}=i_j$. Then (by Wick's lemma and the Cauchy-Schwarz inequality) :
\BEA &&  \Var :{\cal R}{\bf B}^{n,\eta}_{ts}(i_1,\ldots,i_n): =
 \esper \left| :{\cal R}{\bf B}^{n,\eta}_{ts}(i_1,\ldots,i_n): \right|^2 \nonumber\\ && =
\sum_{\sigma'\in\Sigma_{\vec{i}}} \esper\left[ :{\cal R}{\bf B}^{n,\eta}_{ts}(1,\ldots,n):
\   :{\cal R}{\bf B}^{n,\eta}_{ts}(\sigma'(1),\ldots,\sigma'(n)):  \right] \nonumber\\
&& \le |\Sigma_{\vec{i}}| \ .\ \esper |{\cal R}{\bf B}^{n,\eta}(1,\ldots,n)|^2, \label{eq:7:3}\EEA
hence the H\"older and rate estimates of section 3  also hold for \\ $:{\cal R}{\bf B}^{n,\eta}(i_1,\ldots,i_n):$.

One must now prove that the estimates of section 3 hold true  for all possible contractions of
${\cal R}{\bf B}^{n,\eta}(i_1,\ldots,i_n)$. Fixing some non-trivial contraction $(j_1,j_2),\ldots,
(j_{2l-1},j_{2l})$, $l\ge 1$, results in an  expression ${\bf X}^{contr}_{ts}$
belonging to the chaos of order $n-2l$. By necessity, $i_{j_1}=i_{j_2},\ldots,i_{j_{2l-1}}=i_{j_{2l}}$, but it
may well be that there are other index coincidences. The same reasoning as in the case of $:{\cal R}{\bf B}_{ts}^{n,\eta}(i_1,\ldots,i_n):$ (see eq. (\ref{eq:7:3}))  shows that one may actually assume $i_m\not=i_{m'}$ if $m\not=m'$
and $\{m,m'\}\not=\{j_1,j_2\},\ldots,\{j_{2l-1},j_{2l}\}$. Now (as we shall presently prove) the
tree integrals related to the contracted iterated integral ${\bf X}^{contr}_{ts}$ 
may be estimated by considering the tree integrals related to
$\check{\bf X}_{ts}:={\cal R}{\bf B}_{ts}^{n-2l,r}(i_1,\ldots,\check{i_{j_1}},\ldots,\check{i_{j_{2l}}},\ldots,
i_n)$ (which has same law as ${\cal R}{\bf B}^{n-2l,r}_{ts}(1,\ldots,n-2l)$) and (following
the idea introduced in the course of the proof of  Lemma \ref{lem:3:increment-Holder-rate}) increasing
by one the weight $W$ of some other (possibly coinciding) indices $j'_1,\ldots,j'_{2l}\not=
j_1,\ldots, j_{2l}$ -- or, in other words,
'inserting' a factor $|\xi_{j'_1}|^{-2\alpha}\ldots |\xi_{j'_{2l}}|^{-2\alpha}$ in the variance integrals --.
 This amounts in the end
to increasing the H\"older regularity $(n-2l)\alpha^-$ 
 of $\check{\bf X}_{ts}$ by $2l\alpha$, which gives the expected regularity.

Fix  some permutation $\sigma\in\Sigma_n$, and consider the integral over the Fourier
domain
$|\xi_{\sigma(1)}|\le \ldots \le |\xi_{\sigma(n)}|$ as in section 2. Change as before  the order of integration and
the names of the indices so that $dW_{\xi_{\sigma(j)}}(i_j)\to dW_{\xi_j}(i_{\sigma(j)})$;
 for convenience, we shall still index the pairing indices as $(j_1,j_2),\ldots,(j_{2l-1},j_{2l})$.
We may assume that $|j_{2k-1}-j_{2k}|= 1$, $k=1,\ldots,l$ (otherwise $|\xi_m|=|\xi_{j_{2k-1}}|=
|\xi_{j_{2k}}|$ for $j_{2k-1}<m<j_{2k}$ or $j_{2k}<m<j_{2k-1}$, which corresponds to
 a Fourier subdomain of
zero Lebesgue measure). In the sequel, we fix $\sigma\in\Sigma_n$ and $(j,j')=(j_{2k-1},j_{2k})$
for some $k$.

\bigskip

Let  $\tilde{\T}=\tilde{\T}_1\ldots\tilde{\T}_L$
 be a forest appearing in the decomposition of the permutation graph $\T^{\sigma}$ as in subsection 2.4.
Applying  repeatedly Lemma \ref{lem:2:reg} to $\tilde{\T}$ leads to a sum of
 terms obtained from the contraction of $A_{ts}=A_{ts}(1)\ldots A_{ts}(L)$, with
$ A_{ts}(k)= [\del{\cal R}\SkI_{L\tilde{\T}_{k}}]_{ts} 
\prod_j [{\cal R}\SkI_{\T'_{k,j}}]_{s} \left(\left(  \otimes_{v\in V(\tilde{\T}_k)}
dB^{\eta}(i_v)\right)^{\sigma}\right)$, 
where   $L\tilde{\T}_k,\T'_{k,1},\ldots,
\T'_{k,j},\ldots$ are all subtrees appearing in the splitting associated to $A_{ts}(k)$ (see proof
of Lemma \ref{lem:3:boundary-Holder-rate}).

Let $\T$ be one of the above trees, either $L\tilde{\T}_k$ or $\T'_{k,j}$. Reconsider the proof of the H\"older
 estimate or rate of convergence in Lemma  
\ref{lem:3:increment-Holder-rate} or Lemma \ref{lem:3:boundary-Holder-rate}. The integrals
$\left[ \SkI \left( (x_v)_{v\in V(\T)} \to e^{\II \sum_{v\in V(\T)} x_v \xi_v} \right) \right]_u$
appearing in the definition of the regularized skeleton integrals write $i^{-|V(\T)|}
\frac{e^{\II u\sum_{v\in V(\T)} \xi_v}}{
\prod_{v\in V(\T)} (\xi_v+\sum_{w\twoheadrightarrow v} \xi_w)}$ (see Lemma
\ref{lem:2:SkI}). After the contractions, one must sum over Fourier
indices $(\xi_v)_{v\in V(\T)}$ such that $(\xi_v)_{v\in V(\T)}\in \R^{\T}_{reg}$ and $\xi_{j_{2m-1}}=-\xi_{j_{2m}}$ if both
$j_{2m-1},j_{2m}\in V(\T)$.

\begin{figure}[h]
  \centering
   \includegraphics[scale=0.35]{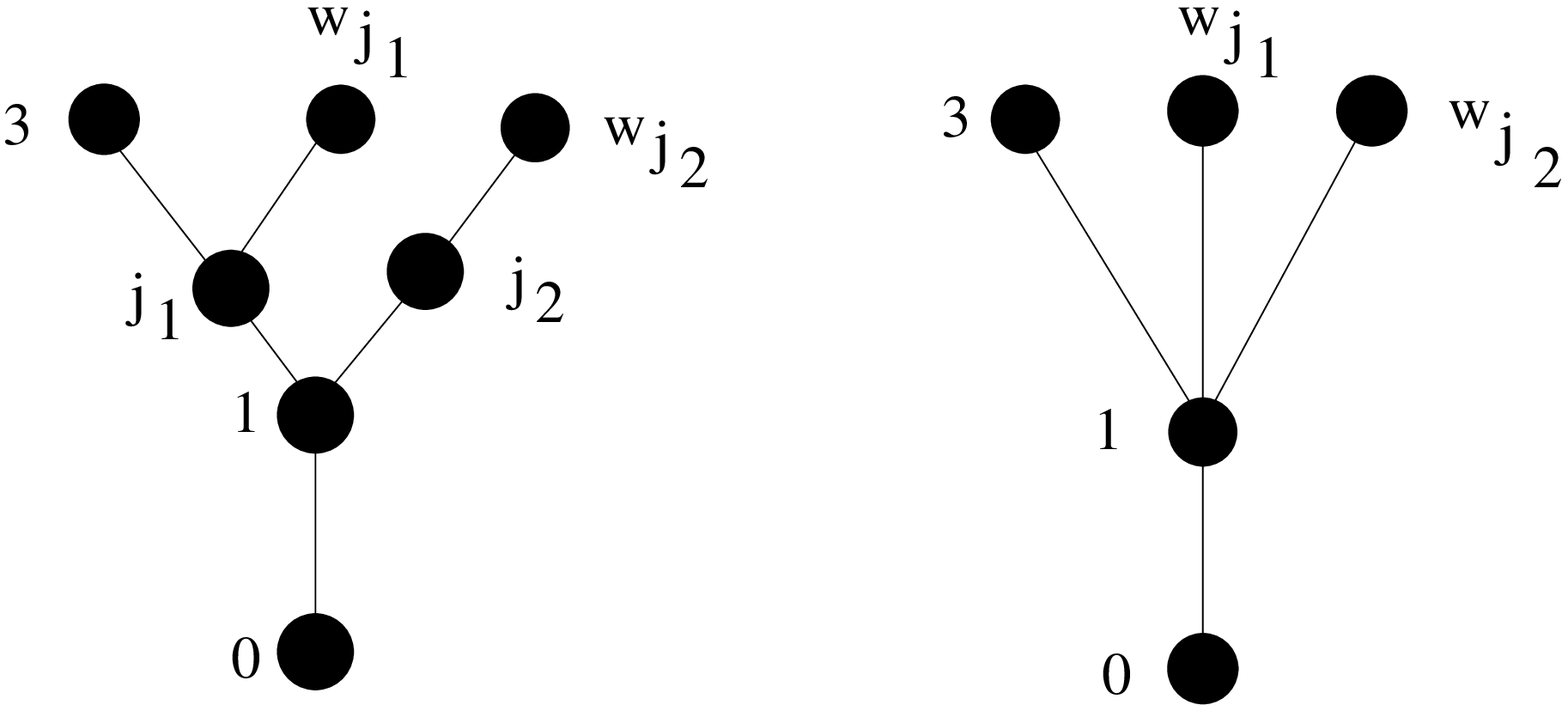}
   \caption{\small{Case (i-a). $\T$ and $\check{\T}$.}}
  \label{Fig7}
\end{figure}

\begin{figure}[h]
  \centering
   \includegraphics[scale=0.35]{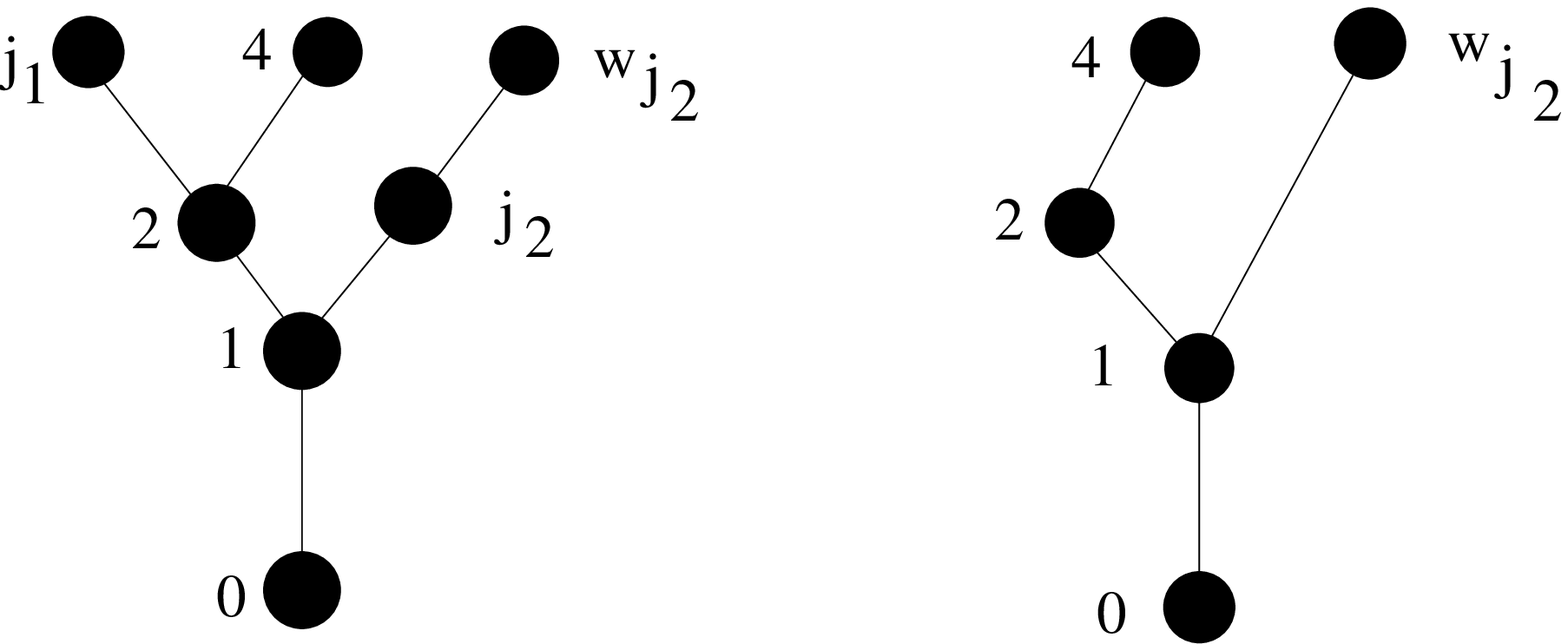}
   \caption{\small{Case (i-b). $\T$ and $\check{\T}$.}}
  \label{Fig8}
\end{figure}

\begin{figure}[h]
  \centering
   \includegraphics[scale=0.35]{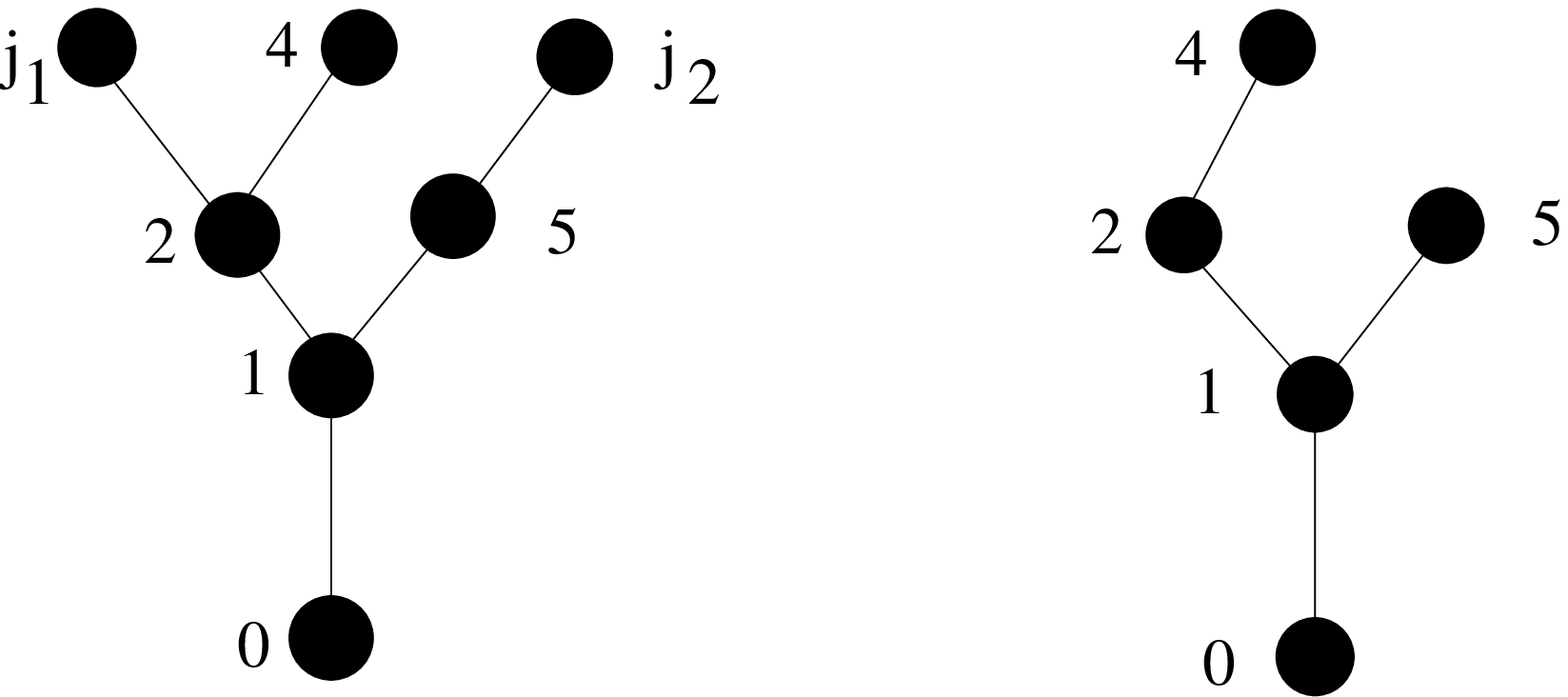}
   \caption{\small{Case (i-c). $\T$ and $\check{\T}$.}}
  \label{Fig9}
\end{figure}

Let \ $\check{\T}$ be the contracted tree
obtained by 'skipping' $\{j_1,\ldots,j_{2l}\}\cap V(\T)$ while going down the tree $\T$
 (see Fig. \ref{Fig7}, \ref{Fig8}, \ref{Fig9}).

The denominator $|\xi_v+\sum_{w\in\T, w\twoheadrightarrow v} \xi_w|$  is
larger (up to a constant) than the denominator  $|\xi_v+\sum_{w\in\check{\T}, w\twoheadrightarrow v} \xi_w|$
obtained by considering the same term in the contracted tree integral $\check{\bf X}_{ts}$ (namely,
$|\xi_v+\sum_{w\in\T, w\twoheadrightarrow v} \xi_w|$ is of the same order as
$\max\{|\xi_w|; w\in\T,w\twoheadrightarrow v\}\ge \max\{|\xi_w|;w\in\check{\T},w\twoheadrightarrow v\}$).
Hence $\esper (A_{ts}^{contr})^2$ may be bounded in the same way as
 $\esper A_{ts}^2$ in the proof of Lemma \ref{lem:3:increment-Holder-rate} or Lemma \ref{lem:3:boundary-Holder-rate}, except that each term in the sum over
$(\xi_{v},v\in V(\T),v\not=j_1,\ldots,j_{2l})$ comes with an extra multiplicative pre-factor
$S=S((\xi_v),v\in V(\T),v\not=j_1,\ldots,j_{2l})$  -- due to the sum over $(\xi_{j_m})_{m=1,\ldots,2l}$ --
which may be seen as an 'insertion'.

Let us estimate this prefactor. We shall assume for the sake of clarity that there is a single contraction
$(j_1,j_2)=(j,j')$
(otherwise the prefactor should be evaluated by contracting each tree  in several stages, 'skipping' successively
$(j_1,j_2),\ldots,(j_{2l-1},j_{2l})$ by pairs). As already mentioned before, $|j-j'|=1$ so that $j$ and $j'$ must be successive
vertices if they belong to the same branch of the same tree $\T$.  Note that, if $j$ and $j'$ are on the same tree,
 the Fourier index 
$\Xi:=\sum_{v\in V(\T)} \xi_v$ (used in the Fourier decomposition of Lemma \ref{lem:3:increment-Holder-rate} or in the
third step of Lemma \ref{lem:3:boundary-Holder-rate}) is left unchanged since $\xi_j+\xi_{j'}=0$.

\medskip

\underline{Case (i)}: $(j,j')$ belong to unconnected branches of the same tree  $\T$.
 This case splits into three
different subcases:

\medskip

\begin{itemize}

\item[(i-a)] neither $j$ nor $j'$ is a leaf. Let $w$, resp. $w'$ be the leaf above  $j$, resp. $j'$ of
maximal index and assume (without loss of generality) that $|\xi_{w}|\le |\xi_{w'}|$. Then 
\BEA && S \lesssim  \left( \int_{|\xi_{j}|\le |\xi_{w}|} d\xi_j \frac{|\xi_j|^{1-2\alpha}}{|\xi_{w}\xi_{w'}|} \right)^2 
\lesssim  \left( \int_{|\xi_{j}|\le |\xi_{w}|}  d\xi_j |\xi_{w}|^{-1-2\alpha} \right)^2 \lesssim |\xi_{w}|^{-4\alpha}\nonumber\\  \EEA
which has the effect of increasing the weight $W(w)$ by $2$.

\item[(i-b)] $j$ is a leaf, $j'$ is not. Let $w'$ be the leaf of maximal index above $j'$. Then
\BEA && S\le \left( \int_{|\xi_{j}|\le |\xi_{w'}|} d\xi_j  \frac{ |\xi_{j}|^{1-2\alpha}}{|\xi_{j} \xi_{w'}|} \right)^2 
\lesssim \left( \frac{1}{|\xi_{w'}|} \int_{|\xi_{j}|\le |\xi_{w'}|} d\xi_j  |\xi_{j}|^{-2\alpha} \right)^2 \lesssim
|\xi_{w'}|^{-4\alpha}. \nonumber\\ 
\EEA

\item[(i-c)] both $j$ and $j'$ are leaves. Let $v$, resp. $v'$  be the vertex below $j$, resp.
$j'$, i.e. $j\to v$, $j'\to v'$.
Then
\BEQ S\lesssim \left( \int_{|\xi_{j}|\ge \max(|\xi_{v}|,|\xi_{v'}|)} d\xi_j |\xi_j|^{-1-2\alpha}
\right)^2\lesssim |\xi_{v}|^{-4\alpha}\EEQ
which has the effect of increasing $W(v)$ by $2$.

\end{itemize}

\underline{Case (ii)}: $(j,j')$ are successive vertices on the same branch of the same tree
$\T$. Assume (without loss of generality)
that $j\to j'$. Then $S=0$ if $j$ is a leaf
 (since $\xi_{j'}+\sum_{w\twoheadrightarrow j'}\xi_w=\xi_{j}+\xi_{j'}=0$
 and such indices fail to meet the condition defining
$\R^{\T}_{reg}$),  otherwise  $S\lesssim |\xi_w|^{-4\alpha}$ if $w$ is the leaf of maximal index above
$j$ (by the same argument as in case (i-a)).

\underline{Case (iii)}: $(j,j')$ belong to two different trees, $\T$ and $\T'$.

This case is a variant of case (i). 
Nothing changes compared to case (i) unless (as in the proof of Lemma \ref{lem:3:increment-Holder-rate}
or in the 3rd step of Lemma \ref{lem:3:boundary-Holder-rate}) one needs to compute the variance
of the coefficient $a(\Xi)$ or $a_s(\Xi)$ of $e^{\II u\Xi}$ for $\Xi$ fixed. Assume $j$  belongs
to the tree $\T=L\tilde{\T}_{k}$ while $j'$ is on one of the cut trees $\T'_{k,1},\ldots,\T'_{k,j},\ldots$

\smallskip

Assume first $j$ is not a leaf, 
 and let $w$  be the leaf above $j$. Then the presence
of the extra vertex $j$ modifies the Fourier index $\Xi$ in the Fourier decomposition of $A^{contr}_{ts}(k)$,
$A^{contr}_{ts}(k)=\int_{\R} a(\Xi) (e^{\II \Xi t}-e^{\II \Xi s}) d\Xi$ or
$A^{contr}_{ts}(k)=\int_{\R} a_s(\Xi) (e^{\II \Xi t}-e^{\II \Xi s}) d\Xi$,
 by a factor which is bounded and bounded away from $0$, hence $S\lesssim |\xi_{w}|^{-4\alpha}$ as in case (i-a).

\smallskip  If $j$ is a leaf as in case (i-b) -- while
$w'$ is as before the leaf of maximal index over $j'$ --, one has: $|\xi_{j}|\lesssim |\Xi|\lesssim |\xi_{j}|$. Hence
the sum over $\xi_{j}$ contributes an extra multiplicative pre-factor $S$ to the variance of
the coefficient of $a(\Xi)$ or $a_s(\Xi)$ of order
\BEQ S\lesssim \left( \int_{|\Xi|/2\le |\xi_j|\le 2|\Xi|} d\xi_j
 \frac{|\xi_{j}|^{1-2\alpha}}{|\xi_{j}\xi_{w'}|}\right)^2 \lesssim
\left(  \int_{|\Xi|/2\le |\xi_j|\le 2|\Xi|}
 |\xi_j|^{-1-2\alpha}  \right)^2 \lesssim |\Xi|^{-4\alpha},\EEQ
which increases the H\"older index by $2\alpha$ (see Lemma \ref{lem:3:Kah}).

\medskip

The case when both $j$ and $j'$ belong to left parts $L\tilde{\T}_{k}$, $L\tilde{\T}_{k'}$
 is similar and left to the
reader. \hfill \eop

\medskip

This concludes at last the proof of Theorem \ref{th:0}.


\subsection{A remark: about the two-dimensional antisymmetric fBm}

Consider a one-dimensional analytic fractional Brownian motion $\Gamma$ as in \cite{TinUnt08}.

\begin{Definition}

{\it Let $Z_t=(Z_t(1),Z_t(2))=(2\Re \Gamma_t,2\Im\Gamma_t)$, $t\in\R$.
We call this new centered Gaussian process indexed by $\R$
 the {\em two-dimensional antisymmetric fBm}.
}

\end{Definition}

Its paths are  a.s. $\alpha^-$-H\"older. The
marginal processes $Z(1)$, $Z(2)$ are usual
fractional Brownian motions.  The covariance between
$Z(1)$ and $Z(2)$ writes
(see \cite{TinUnt08}) 
\BEQ \Cov(Z_s(1),Z_t(2))=-\frac{\tan\pi\alpha}{2}[-\sgn(s)|s|^{2\alpha}+\sgn(t) |t|^{2\alpha}
-\sgn(t-s) |t-s|^{2\alpha}].\EEQ

Note that we never used any particular linear combination of the analytic/anti-analytic components of
$B$ in the estimates of section 3 and 4. Hence these also hold for  $Z$, which gives for
free a rough path over $Z$  satisfying Theorem \ref{th:0} of the Introduction.


\end{document}